\documentclass[final,reqno]{amsart}

\usepackage{graphicx}
\usepackage{amsmath}

\usepackage{epstopdf} 
\usepackage{pdfsync}
\usepackage{amsmath,amssymb}
\usepackage{enumerate}
\usepackage{xspace}
\usepackage{boxedminipage}
\usepackage{algorithmicx}
\usepackage[ruled]{algorithm}
\usepackage{algpseudocode}
\usepackage{listings}
\usepackage{array}
\usepackage{tabularx}
\usepackage{subfigure}
\usepackage{caption}
\usepackage[cm]{fullpage}
\usepackage{james}

\usepackage{tikz}
\usetikzlibrary{calc}
\usepackage{pgfplots}
\pgfplotsset{width=7cm,compat=1.15}
\usepgfplotslibrary{colorbrewer}

\author{
  James Jackaman
}
\address{
  James Jackaman
  \thanks{
    Department of Mathematics and Statistics, Memorial University of
    Newfoundland, St.\ John's, NL, A1C 5S7, Canada
    {\tt{jjackaman@mun.ca}}}
}

\author{
  Tristan Pryer
}
\address{
  Tristan Pryer
  Department of Mathematical Sciences,
  University of Bath, Bath BA2 7AY, UK
  {\tt{tmp38@bath.ac.uk}}
}

\thanks{ J.J. was partially supported through a PhD scholarship
  awarded by the ``EPSRC Centre for Doctoral Training in the
  Mathematics of Planet Earth at Imperial College London and the
  University of Reading'' EP/L016613/1, and the Canadian Research
  Chairs and NSERC Discovery grant programs. T.P. was partially
  supported through the EPSRC grant EP/P000835/1 and the Newton Fund
  grant 261865400. Both authors would additionally like to acknowledge
  the support of the Isaac Newton Institute for Mathematical Sciences,
  Cambridge through the EPSRC grant EP/K032208/1.  }

\title[]{Conservative Galerkin methods for dispersive Hamiltonian
  problems}

\begin{document}

\maketitle

\begin{abstract}

  An energy conservative discontinuous Galerkin scheme for a
  generalised third order KdV type equation is designed. Based on the
  conservation principle, we propose techniques that allow for the
  derivation of optimal a priori bounds for the linear KdV equation
  and a posteriori bounds for the linear and modified KdV
  equation. Extensive numerical experiments showcasing the good long
  time behaviour of the scheme are summarised which are in agreement
  with the analysis proposed.

\end{abstract}

\section{Introduction}
\label{sec:introduction}

Hamiltonian partial differential equations (PDEs) arise naturally from
a variety of physically motivated application areas, with decisive
examples including meteorological, as exemplified by the semi-geostrophic equations
\cite{RoulstoneNorbury:1994}, and oceanographical, such as the
Korteweg-de Vries (KdV) and nonlinear Schr\"odinger equations
\cite{MullerGarretOsborne:2005}. The KdV and nonlinear Schr\"odinger
equations are particularly special examples, in that they are
bi-Hamiltonian \cite{Magri:1978}. This means they have two different
Hamiltonian formulations which, in turn, is one way to understand the
notion of \emph{integrability} of these problems.

Conservative schemes for Hamiltonian ordinary differential equations
(ODEs) are, by now, relatively well understood, see
\cite[c.f.]{Celledoni:2012, LeimkuhlerReich:2004,
  HairerLubichWanner:2006, BokhoveLynch:2007,
  BlanesCasa:2016}. Typically numerical schemes designed for this
class of problem have some property of the ODE built into them, for
example preservation of the Hamiltonian or the underlying symplectic
form, are classified as geometric integrators.

In this contribution, we consider a family of Hamiltonian PDEs that
generalise the famous defocusing modified KdV equation
\begin{equation*}
  u_t
  -
  \frac 3 2 u^2 u_x
  +
  u_{xxx} = 0,
  \label{eq:mkdv}
\end{equation*}
where the sub-indices denote partial differentiation with respect to
the corresponding independent variable. This equation has numerous
applications including fluid dynamics and plasma physics
\cite{AblowitzClarkson:1991}. The Korteweg third order
term, as well as the non-linearity, can cause issues in the numerical
treatment of this problem. In particular, the combination of these two
terms cause significant trouble in the design of numerical
methods that are optimally convergent in the function approximation
sense.

In previous numerical studies of the scalar KdV and modified KdV
equations \cite[c.f.]{XuShu:2007,YanShu:2002}, it has been observed
that classical finite volume and discontinuous Galerkin (dG) schemes
with ``standard'' numerical fluxes introduce numerical
artefacts. Consequently the long-term dynamics of solutions may be
destroyed by the addition of \emph{artificial diffusion}. Such
diffusion in a given scheme endows it with desirable stability
properties, however, it typically destroys all information in the
long-term dynamics of the system through smearing of solutions.

Hamiltonian problems are inherently conservative in the sense that the
underlying Hamiltonian is conserved over time. Such equations may have
additional structures which manifest themselves through additional
conserved quantities. In particular, for the KdV equation mass and
momentum are such quantities. In
\cite{BonaChenKarakashianXing:2013,KarakashianMakridakis:2015} the
authors propose and analyse a dG method for generalised KdV
equations. The method itself is very carefully designed to be
conservative, in that the invariant corresponding to the
\emph{momentum} is inherited by the discretisation. This naturally
yields $\leb{2}$ stability for the numerical method along with
extremely good long time dynamics. In the scalar case one may also
design schemes that conserve the energy itself
\cite{Winther:1980,JackamanPapamikosPryer:2017}, however, it does not
seem possible to design schemes to conserve more than two of these
invariants. \revise{This is, in part, due to the invariants differing
  in order of nonlinearity}. Other promising mechanisms to approximate
such problems include those based on the Fokas transform
\cite{KesiciPelloniPryerSmith:2018} and the Chebfun package
\cite{driscoll2008chebop} both of which experimentally has extremely
good long time properties as the schemes naturally inherit many of the
properties of the PDE up to small precision. \revise{Alternatively,
  one may aim to preserve the multi-symplectic structure of the KdV
  equation, which has proven successful with Preissmann box schemes
  \cite{AscherMcLachlan:2005, ZhaoQin:2000} resulting in
  discretisations possessing desirable qualitative behaviour. The
  study of multi-symplectic schemes is a fertile area of research with
  finite element discretisations utilising this structure currently
  being developed, for example \cite{McLachlanStern:2020,
    self:multisym}.  }

Our goal in this work is the derivation of Galerkin discretisations
aimed at preserving the underlying algebraic properties satisfied by
the PDE system whilst avoiding the introduction of stabilising
diffusion terms. Our schemes are therefore consistent with (one of)
the Hamiltonian formulation(s) of the original problem which
physically represents energy. We note that our approach is not an
adaptation of entropy conserving schemes developed for systems of
conservation laws, rather we study the algebraic properties of the PDE
and formulate the discretisation to inherit this specific
structure. Our methods are of arbitrarily high order accuracy in
space, extendable to arbitrarily high order in time, and provide
relevant approximations free from numerical artefacts. Similar
techniques have proven useful in the study of dispersive phase flow
problems \cite{GiesselmannPryer:2015a,
  GiesselmannMakridakisPryer:2014} and we anticipate they will be
extremely useful in dynamic model adaptivity
\cite{GiesselmannPryer:2017a}.

To highlight the good behaviour of the scheme we propose, we give an a
priori error analysis for the linear problem showing optimal error
bounds in the energy norm, under the assumption of smoothness in the
initial conditions. Further, we give an a posteriori error analysis
making use of a hybrid framework consisting of elliptic reconstruction
techniques \cite{MN06,Lakkis:2006, LakkisMakridakisPryer:2015}
together with those developed for hyperbolic conservation laws
\cite{GiesselmannMakridakisPryer:2015} to allow us to derive optimal a
posteriori error bounds in the energy norm. Note that the arguments we
use are quite different to that of \cite{KarakashianMakridakis:2015}
where the authors construct a dispersive reconstruction to allow for a
posteriori control in $\leb{2}$.

The remainder of this work is set out as follows: In Section \ref{sec:setup}
we introduce notation, the model problem and some of its
properties. We also give some exact solutions to the problem that will
be useful in numerical benchmarking. In Section \ref{sec:dis} we reformulate
of the problem through the introduction of an auxiliary variable,
introduced to allow for a simple construction of the numerical
scheme. We propose a spatial discretisation based on discontinuous
Galerkin finite elements, show it is conservative, well posed, and
give an a priori analysis yielding optimal error bounds in the energy
norm. In Section \ref{sec:apos} we give an a posteriori analysis making use
of elliptic reconstruction techniques. In
Section \ref{sec:temp-and-numerics} we state a fully discrete scheme, show
it is conservative and summarise extensive numerical experiments
validating the analysis done in Section \ref{sec:dis}--\ref{sec:apos}.

\section{Problem setup}
\label{sec:setup}

In this section we formulate the model problem, fix notation and give
some basic assumptions. We describe some known results and history of
the defocusing generalised Korteweg-de Vries equation, highlighting
the Hamiltonian structure of the equation. We show that the underlying
Hamiltonian structure naturally yields an induced stability of the
solutions to the PDE system and give a summary of some exact solutions
for specific non-linearities.

Throughout this work we consider the (1+1)-dimensional dispersive
Cauchy problem
\begin{equation}
  \label{thorder}
  \begin{split}
    u_t - f'(u)_x + u_{xxx} = 0& \qquad x\in S^1, \qquad t>0
    \\
    u(x,0) = u^0(x)& \qquad x\in S^1
  \end{split}
\end{equation}
with periodic boundary conditions over the unit interval $S^1:=[0,1]$
with endpoints being identified with one another. The non-linearity is
polynomial, specifically
\begin{equation}
  f(u) = \alpha u^m  \text{ for } 2 \le m\in \mathbb{Z}, 0 < \alpha \in \real.
\end{equation}
Examples of the PDE include the linear KdV equation
\begin{equation}
  u_t - \alpha u_x + u_{xxx} = 0,
\end{equation}
and the mKdV equation
\begin{equation}
  u_t - \alpha u^2 u_x + u_{xxx} = 0.
\end{equation}
These equations are the focus of the analysis in this work.

Notice the sign in front of the first order term. In an analogy to the
nonlinear Schr\"odinger equation we refer to as a defocusing equation,
with the focusing equations having the opposing sign on the
non-linearity.

\begin{Proposition}
  The dispersive problem (\ref{thorder}) has the following three base
  invariants:
  \begin{equation} \label{eqn:invariants}
    \begin{split}
      \ddt F_{-1}(u) &:= \ddt \int_{S^1} u \d x = 0
      \\
      \ddt F_0(u)
      &:=
      \ddt \int_{S^1} \frac 1 2 u^2 \d x = 0
      \\
      \ddt F_1(u)
      &:=
      \ddt \int_{S^1} \frac 1 2 u_x^2 + f(u)  \d x = 0.
    \end{split}
  \end{equation}
\end{Proposition}

Throughout this work we denote the standard Lebesgue spaces by
$\leb p(\omega)$, $1\le p\le \infty$, $\omega\subset\mathbb{R}$, with
corresponding norms $\|\cdot\|_{L^p(\omega)}$. Let also $H^s(\omega)$,
be the Hilbertian Sobolev space of index $s\in\mathbb{R}$ of
real-valued functions defined on $\omega\subset\mathbb{R}^d$,
constructed via standard interpolation and/or duality procedures,
along with the corresponding norm and semi-norm
\begin{gather}
  \Norm{u}_{\sob{k}{p}(\W)}
  := 
  \begin{cases}
    \qp{\sum_{\norm{\alpha}\leq k}\Norm{\D^{\alpha} u}_{\leb{p}(\W)}^p}^{1/p} &\text{ if } p \in [1,\infty)
    \\
    \sum_{\norm{\alpha}\leq k}\Norm{\D^{\alpha} u}_{\leb{\infty}(\W)} &\text{ if } p = \infty 
  \end{cases}
  \\
  \norm{u}_{\sob{k}{p}(\W)}
  :=
  \Norm{\D^k u}_{\leb{p}(\W)} 
\end{gather}
respectively. We also make use of the following notation for time
dependent Sobolev (Bochner) spaces:
\begin{gather}
  \cont{i}(0,T; \sobh{k}(S^1))
  :=
  \ensemble{u : [0,T] \to \sobh{k}(S^1)}
           {u \text{ and temporal derivatives up to $i$-th order are continuous}},
           \\
           \leb{\infty}(0,T; \sobh{k}(\W))
  :=
           \ensemble{u : [0,T] \to \sobh{k}(\W)}
                    {\esssup_{t\in [0,T]} \Norm{u(t)}_{\sobh{k}(\W)} < \infty}.
\end{gather}

Under some regularity assumptions on the initial condition one can
make use of semi-group techniques to show the following:

\begin{Proposition}[(\ref{thorder}) is well posed \cite{Kato:1975}]
  Given $u^0 \in \sobh s(S^1)$, with $s\geq 3$. Then (\ref{thorder})
  has a unique solution with
  \begin{equation}
    u \in \cont{0}(0,T; \sobh s(S^1)) \cap \cont{1}(0,T; \sobh{s-3}(S^1)),
  \end{equation}
  for arbitrary $T>0$.
\end{Proposition}

\begin{Proposition}[Pointwise solution control]
  \label{pro:pointwise}
  Notice that energy conservation
  \begin{equation} 
    \ddt \int_{S^1} u_x^2 + f(u) \d x = 0 
  \end{equation}
  and mass conservation
  \begin{equation}
    \ddt \int_{S^1} u \d x = 0
  \end{equation}
  immediately shows that,
  \begin{equation}
    \Norm{u}_{\cont{i}(0,T; \sobh{1}(S^1))}
    \leq    
    C\qp{\Norm{u^0_x}_{\leb 2(S^1)}^2
    +
    \Norm{u^0}_{\leb{m}(S^1)}^{m}}^{1/2},
  \end{equation}
  for some constant $C>0$.
  Since $\sobh1(S^1) \subset \leb{\infty}(S^1)$ we see
  \begin{equation}
    \sup_{t\in [0,\infty]} \Norm{u(t)}_{\leb{\infty}(S^1)}
    \leq
    C\qp{
    \Norm{u^0_x}_{\leb 2(S^1)}^2
    +
    \Norm{u^0}_{\leb{m}(S^1)}^{m}}^{1/2}.
  \end{equation}
\end{Proposition}

\begin{Proposition}[Exact solution to the linear problem]
  Let $f(u) = \frac12 u^2$, then under the ansatz that $u(t,x) =
  u(\xi)$, with $\xi = c \left(x + (1+c^2)t\right)$ we find that
  \begin{equation} \label{eqn:ldefocusexact}
    u(x,t) 
    =
    C_1 \sin{\xi}
    + C_2 \cos{\xi}
    ,
  \end{equation}
  solves \eqref{thorder} where $c = 2 l \pi $ for $l \in \mathbb{Z}$
  and $C_1,C_2$ denote real constants. Due to the linear nature of the
  problem any linear combination of \eqref{eqn:ldefocusexact} for
  various attainable parameter values is also a solution.
\end{Proposition}

\begin{Proposition}[Exact solution to the nonlinear problem]
  With $f(u) = 6 u^3$, then it can be shown that the position solution
  \begin{equation}
    u(x,t) =
    \frac c 2 \operatorname{csch}\qp{-\frac {c^{1/2}\qp{x-ct}} 2}^2
  \end{equation}
  formally solves (\ref{thorder}). It is well-known that one can map
  solutions from the defocusing mKdV equation to solutions to the KdV
  equation employing the Miura transform. Although, it is worth noting
  that it is not possible to get smooth, non-singular position solutions
  of the defocusing mKdV through inverse scattering techniques because
  of the singularity that is inherent in its Darboux
  transformation. For $f(u) = \frac 1 4 u^4$, one can, however find
  kink
  \begin{equation}
    u(x, t) =
    (3c)^{1/2}\tanh\qp{\frac{(2c)^{1/2}} 2 (x+ct)}
  \end{equation}
  and anti-kink solutions
  \begin{equation}
    u(x, t) =
    -(3c)^{1/2}\tanh\qp{\frac{(2c)^{1/2}} 2 (x+ct)},
  \end{equation}
  that are smooth, but are not periodic. To establish periodic, smooth
  exact solutions, one must examine Jacobi elliptic functions
  \cite{Pava:2006}. Let $sn(x,k)$ denote that Jacobi elliptic function
  with modulus $k \in [0,1)$, then, with $f(u) = \frac 1 2 u^4$, a
    solution is given by \cite{Deconinck:2011}
  \begin{equation}
    \label{eqn:defocusexact}
    u(x,t) = k \operatorname{sn}(x+(k^2+1)t,k).
  \end{equation}
\end{Proposition}

\section{Discretisation and a priori analysis}
\label{sec:dis}

We consider the approximation of \eqref{thorder} by a semi-discrete
discontinuous Galerkin scheme. Let $0 = x_0 < x_1 < \dots < x_N = 1$
be a partition of the periodic domain $S^1$. We denote
$I_j=[x_j,x_{j+1}]$ to be the $j$--th sub-interval and let $h_j:=
x_{j+1}-x_j$ be its size. We denote the piece-wise constant mesh-size
function $h : S^1 \rightarrow [0,\infty)$ where $h |_{I_j}=h_j$. For
  the purposes of this work, we will assume that $\max{\left( h_j
    N\right) \leq C}$ for some $C>0$. For $q \geq 1$ let $\poly q(I)$
  be the space of polynomials of degree less than or equal to $q$ {on
    $I$}, then we denote
\begin{equation}
  \fes_q
  :=
  \ensemble{g : S^1 \to \rR }{ g \vert_{I_j}
     \in \poly{q}{(I_j)} \text{ for }  j =0,\dots, N-1
  }.
\end{equation}
In addition, we define jump and average operators by
\begin{equation}
  \begin{split}
    \jump{g}_j
    &:= 
    g(x_j^-) - g(x_j^+)
    := 
    \lim_{s \searrow 0} g(x_j-s) - \lim_{s \searrow 0}  g(x_j+s),
    \\
    \avg{g}_j
    &:
    = 
    \frac{1}{2} \qp{ g(x_j^-) +g(x_j^+)}
    :=
    \frac{1}{2} \qp{\lim_{s \searrow 0} g(x_j-s) + \lim_{s \searrow 0} g(x_j+s)}
  \end{split}
\end{equation}
where the periodic boundary conditions are accounted for by
$\jump{g}_0 := g(x_N^-)-g(x_0^+)$ and $\avg{g}_0 := \frac 1 2
\qp{g(x_N^-)+g(x_0^+)}.$ Throughout this work we will use the
convention that $C > 0$ denotes a generic constant which may depend on
$q$, the ratio of concurrent cell sizes and non-linearity degree $m$,
but is independent of $h$ and the exact solution $u$.

We will examine semi-discrete numerical schemes which are based on the
following reformulation of \eqref{thorder} using an auxiliary variable
$v$
\begin{equation}
  \label{eq:mixed}
  \begin{split}
    u_t +v_x &= 0
    \\
    v + f'(u) - u_{xx} &= 0.
  \end{split}
\end{equation}
The purpose of this variable becomes apparent in the discretisation of
(\ref{thorder}). Indeed, $v$ is deliberately chosen as the \emph{first
  variation} of the energy functional, $F_1$.  We note that a similar
numerical procedure was applied to a regularised elastodynamics
problem in \cite{GiesselmannPryer:2015,GiesselmannPryer:2016}. We will
begin by introducing some projection operators and describe some of
their properties we will make use of throughout this work.

\begin{Definition}[$\leb{2}$ projection operator and properties]
  We define the $\leb{2}$ projection operator
  $\cP_h:L^2(S^1)\to\fes_q$ by requiring
  \begin{equation}
    \sum_{j=0}^{N-1}\int_{x_{j}}^{x_{j+1}} \cP_h(w) \Phi \d x
    =
    \sum_{j=0}^{N-1}\int_{x_{j}}^{x_{j+1}} w \Phi \d x \Foreach \Phi\in\fes_q.
  \end{equation}
  When $w\in\sobh{q+1}(S^1)$, the following approximation
  properties hold
  \begin{equation}
    \Norm{w - \cP_h(w)}_{\leb{2}(S^1)} \leq C h^{q+1} \norm{w}_{\sobh{q+1}(S^1)}.
  \end{equation}
\end{Definition}

\begin{Definition}[Discrete gradients $\cG_h$ and properties]
  The discrete gradient operator $\cG_h: \prod_{j=0}^{N-1} \sobh1(I_j)
  \to \fes_q$ is defined by
  \begin{equation}
    \label{dgrad}
    \sum_{j=0}^{N-1}\int_{x_{j}}^{x_{j+1}} \cG_h(w)\Psi \d x 
    =
    \sum_{j=0}^{N-1} \int_{x_j}^{x_{j+1}} w_x \Psi \d x 
    -
    \sum_{j=0}^{N-1} \jump{w}_j \avg{\Psi}_j
    \quad 
    \forall \ \Psi \in \rV_q.
  \end{equation}
  It can be seen from the definition that these operators satisfy a
  discrete integration by parts, that is for $W,\Psi\in\fes_q$
  \begin{equation}
    \sum_{j=0}^{N-1} \int_{x_j}^{x_{j+1}} \cG_h(W)\Psi \d x 
    =
    - \sum_{j=0}^{N-1} \int_{x_j}^{x_{j+1}} W \cG_h(\Psi) \d x.
  \end{equation}
\end{Definition}

\begin{Definition}[Interior penalty bilinear form and properties]
  We define the interior penalty bi-linear form for $w,\psi \in
  \prod_{j=0}^{N-1} \sobh2(I_j)$ as
  \begin{equation}
    \label{def:ip}
    \bih{w}{\psi}
    := 
    \sum_{j=0}^{N-1} \Big( \int_{x_j}^{x_{j+1}} w_x \psi_x \d x 
    -
    \jump{w}_{j} \avg{\psi_x}_{j}
    -
    \jump{\psi}_{j} \avg{w_x}_{j} + \sigma \avg{h}_j^{-1} \jump{w}_{j} \jump{\psi}_{j}\Big),
  \end{equation}
  for some $\sigma \gg 1$. Note that this is symmetric, that is
  \begin{equation}
    \bih{w}{\psi} = \bih{\psi}{w},
  \end{equation}
  and a consistent representation of the Laplacian so for
  $u\in\sobh{2}(S^1)$ we have
  \begin{equation}
    \bih{u}{\psi} = \int_{S^1} -u_{xx} \psi \d x.
  \end{equation}
\end{Definition}

\subsection*{\bf Semi discrete scheme}
With these definitions in hand we are now in a position to state the
semi discretisation of (\ref{eq:mixed}). This is to seek $U \in
\cont{1}([0,T),\rV_q)$ and $ V \in \cont{0}([0,T),\rV_q)$ such that
\begin{equation}\label{sds}
  \begin{split}
    \int_{S^1} U_t \Phi + \cG_h(V)\Phi \d x 
    &=0 
    \quad \forall \ \Phi \in \rV_q, \ t\in [0,T]\\
    \int_{S^1} \qp{V + f'(U)} \Psi + \bih{U}{\Psi} 
    &=0
    \quad \forall \ \Psi \in \rV_q, \ t\in [0,T]\\
    U(0) &= \cP_h(u^0).
  \end{split}
\end{equation}

\begin{Proposition}[Conservativity of discrete invariants]
  \label{pro:cons}
  Solutions $U \in \cont{1}([0,T),\rV_q)$ and $ V \in
    \cont{0}([0,T),\rV_q)$ to the discrete scheme (\ref{sds}) satisfy
      conservation of mass,
      \begin{equation}
        \label{eq:dismass}
        \ddt F_{-1}(U) = \ddt \qp{\int_{S^1} U \d x} = 0 
      \end{equation}
      and the discrete energy identity
      \begin{equation}
        \label{eq:disenergy}
        \ddt F_{1,h}(U) := \ddt \qp{\frac 12 \bih{U}{U} + \int_{S^1} f(U) \d x} = 0.
      \end{equation}
\end{Proposition}
\begin{proof}
  We see (\ref{eq:dismass}) by taking $\Phi = 1$ in (\ref{sds}) and
  utilising the definition of $\cG_h$. For (\ref{eq:disenergy}) we
  explicitly compute the time derivative and use the symmetry of
  $\bih{\cdot}{\cdot}$. Hence
  \begin{equation}
    \begin{split}
      \ddt F_{1,h}(U)
      &=
      \bih{U}{U_t}
      +
      \int_{S^1} f'(U) U_t \d x
      =
      - \int_{S^1} V U_t \d x
      =
      \int_{S^1} \cG_h(V) V \d x = 0,
    \end{split}
  \end{equation}
  as required.
\end{proof}

\begin{Corollary}[Pointwise discrete solution control]
  Through similar arguments as Proposition \ref{pro:pointwise} we have
  that
  \begin{equation}
    \sup_{t\in [0,T]} \Norm{U(t)}_{\leb{\infty}(S^1)}
    \leq
    C.
  \end{equation}
\end{Corollary}

\begin{Lemma}[Existence and uniqueness to the discrete scheme (\ref{sds})]
  For given initial data $U(0) \in \fes_q$ the ODE system \eqref{sds}
  has a unique solution with $U \in \cont{1}((0,T),\fes_q)).$
\end{Lemma}
\begin{proof}
  We begin by eliminating the auxiliary variable by writing
  (\ref{sds}) in primal form. To that end, we define the discrete
  Laplacian $A_h : \fes_q \to \fes_q$ such that for any fixed
  $\Psi\in\fes_q$
  \begin{equation}
    \int_{S^1} - A_h(\Psi) \Phi \d x = \bih{\Psi}{\Phi} \Foreach \Phi\in\fes_q.
  \end{equation}
  Then (\ref{sds}) can be written as
  \begin{equation}
    \int_{S^1} \qp{U_t - \cG_h(\cP_h (f'(U)) - A_h(U)) }\Phi = 0,
  \end{equation}
  which allows us to interpret the scheme point-wise as an ODE
  \begin{equation}
    U_t =: y'(t) = F(y(t)) := \cG_h(\cP_h (f'(U)) - A_h(U)).
  \end{equation}
  In view of inverse estimates and the stability of the $\leb{2}$
  projector we see that $F$ is continuous. Further, through the
  conservativity of the scheme from Proposition \ref{pro:cons} we see
  that $y$ remains in a bounded set, which depends upon the initial
  data, as long as a classical solution to (\ref{thorder}) exists,
  irrespective of the non-linearity. That is, for $y(0) \in K \subset
  \fes_q$, $y(t) \in K$ for all $t$. Further the Jacobian
  \begin{equation}
    \D F(y)(z) = \cG_h(\cP_h (f''(y)z) - A_h(z))
  \end{equation}
  is a uniformly bounded operator. We
  may now invoke the Picard-Lindel\"of theorem yielding a global solution.
\end{proof}

\subsection*{\bf A priori error analysis}
We dedicate the rest of this section to the a priori error analysis of
the scheme (\ref{sds}) for the linear problem. We proceed by making
use of the discrete stability framework satisfied by the
approximation, introducing appropriate projection operators and
defining the mesh dependent norms for our analysis.

\begin{Lemma}[Perturbed error equation]
  \label{lem:drrer}
  Let $(U,V)$ be a solution of \eqref{sds} and let
  \begin{equation}
    (\widetilde U,\widetilde V) \in \cont{1}([0,T),\rV_q)
      \times \cont{0}([0,T),\rV_q)
  \end{equation}  satisfy the following perturbed problem
  \begin{equation}
    \label{sds:pert}
    \begin{split}
      \int_{S^1} \widetilde U_t \Phi
      + \cG_h(\widetilde V)\Phi \d x
      &=
      \int_{S^1} -\Ru \Phi \d x \quad \forall \ \Phi \in \rV_q \\
      \int_{S^1} \widetilde V \Psi
      + f'(U) - f'(\theta^u) \Psi \d x
      + \bi{\widetilde U}{\Psi}
      &=\int_{S^1} -\Rv \Psi \d x \quad \forall \ \Psi \in \rV_q
      ,
    \end{split}
  \end{equation}
  where $\Ru,\Rv \in \cont{0}([0,T),\rV_q)$ represent discrete
  residuals and $\theta^u = U-\widetilde U$. Then, with
  $\theta^v = V-\widetilde V$ we have
    \begin{equation}
      \ddt F_{1,h}(\theta^u)
      =
      \int_{S^1} f'(\theta^u)  \Ru - \Rv \cG_h(\theta^v) \d x
      +
      \bih{\theta^u}{\Ru}
      .
    \end{equation}
\end{Lemma}

\begin{proof}
  To begin we note that a discrete error equation is given by taking
  the difference of~(\ref{sds}) and (\ref{sds:pert}) yielding
  \begin{equation}
    \label{eq:err-eq}
    \begin{split}
      \int_{S^1} \theta^u_t \Phi + \cG_h(\theta^v) \Phi \d x &= \int_{S^1} \Ru \Phi \d x
      \\
      \int_{S^1} \theta^v \Psi + f'(\theta^u)\Psi \d x + \bih{\theta^u}{\Psi} &= \int_{S^1} \Rv \Psi \d x.
    \end{split}
  \end{equation}
  Explicitly computing the time derivative
  \begin{equation}
    \begin{split}
      \ddt F_{1,h}(\theta^u)
      =
      \ddt \qp{\frac 12 \bih{\theta^u}{\theta^u}
        +
        \int_{S^1} f(\theta^u) \d x}
      &=
      \bih{\theta^u}{\theta^u_t}
      +
      \int_{S^1} f'(\theta^u) \theta^u_t \d x
      .
    \end{split}
  \end{equation}
  Now making use of (\ref{eq:err-eq}) with $\Psi = \theta_t^u$ we see
  \begin{equation}
    \ddt F_{1,h}(\theta^u) = \int_{S^1} \qp{\Rv - \theta^v} \theta^u_t \d x
    .
  \end{equation}
  Again using (\ref{eq:err-eq}), this time with
  $\Phi = \Rv - \theta^v$, we see
  \begin{equation}
    \begin{split}
      \ddt F_{1,h}(\theta^u)
      &=
      \int_{S^1} \qp{\Rv - \theta^v} \qp{\Ru - \cG_h(\theta^v)} \d x
      \\
      &=
      \int_{S^1} \Rv\Ru - \theta^v \Ru - \Rv \cG_h(\theta^v) \d x,
    \end{split}
  \end{equation}
  where we have used skew-symmetry of $\cG_h$. Further, again by
  (\ref{eq:err-eq}) with $\Psi = \Ru$, we have
  \begin{equation}
    \begin{split}
      \ddt F_{1,h}(\theta^u)
      &=
      \int_{S^1} f'(\theta^u)  \Ru - \Rv \cG_h(\theta^v) \d x
      +
      \bih{\theta^u}{\Ru},
    \end{split}
  \end{equation}
  as required.
\end{proof}

\begin{Lemma}[Projection operator $\cS_h$ and error control]
  \label{lem:S}
  Suppose $v\in\sobh{q+2}(S^1)$. Let the polynomial degree $q$ be
  even, the mesh-size be uniform, and the number of elements in the
  mesh, $N$, be odd.  Then, there exists a uniquely defined projection
  operator, $\cS_h(v)\in\fes_q$, satisfying
  \begin{equation}
    \begin{split}
      \int_{S^1} \cS_h(v) \Phi \d x
      &=
      \int_{S^1} v \Phi \d x \Foreach \Phi \in\fes_{q-1}
      \\
      \avg{\cS_h(v)}_j
      &=
      v(x_j) \Foreach j \in [0,N].
    \end{split}
  \end{equation}
  Furthermore, the following error bound holds:
  \begin{equation}
    \label{eq:err-v}
    \Norm{v - \cS_h(v)}_{\leb{2}(S^1)}
    +
    \Norm{v_x - \cG_h(\cS_h(v))}_{\leb{2}(S^1)}
    \leq
    C h^{q+1} \Norm{v}_{\sob{q+2}{\infty}(S^1)}.
  \end{equation}  
\end{Lemma}
\begin{proof}
  We begin by introducing a related projection. Let
  $\cT_h(v)\in\fes_q$ be defined by
  \begin{equation}
    \begin{split}
      \int_{S^1} \cT_h(v) \Phi \d x
      &=
      \int_{S^1} v \Phi \d x \Foreach \Phi \in\fes_{q-1}
      \\
      \cT_h(v)(x_j^+)
      &=
      v(x_j) \Foreach j \in [0,N].
    \end{split}
  \end{equation}
  Notice that, in contrast to $\cS_h$, the projector $\cT_h$ has
  ``one-sided'' boundary conditions, which means it is locally
  constructed. It is uniquely defined and has the approximation
  property
  \begin{equation}
    \label{eq:T-approx}
    \Norm{v - \cT_h(v)}_{\leb{\infty}(I_j)}
    \leq
    C h^{q+1} \Norm{v}_{\sob{q+1}{\infty}(I_j)}.
  \end{equation}
  Proofs of this can be found in \cite[Lem 8]{GiesselmannPryer:2016}.

  The remainder of the proof takes inspiration from \cite[Prop
    3.1]{BonaChenKarakashianXing:2013}. To show properties for $\cS_h$
  we will consider the error $e = \cS_h(v) - \cT_h(v)$. Notice that
  this satisfies the error relations
  \begin{equation}
    \begin{split}
      \int_{x_j}^{x_{j+1}} e \Phi \d x&= 0 \Foreach \Phi \in \poly{q-1}(I_j), j=0, \dots, N-1
      \\
      \avg{e}_j &= \frac12 \qp{v(x_j) - \cT_h(x_j^-)} \text{ for } j=0, \dots, N-1.
    \end{split}
  \end{equation}
  
  Let $l_k\in\poly{k}(-1,1)$ denote the $k^{th}$ Legendre polynomial
  on $(-1,1)$ and $l_{j,k}\in\poly{k}(I_j)$ the transformation to $I_j$ given by
  \begin{equation}
    l_{j,k}(x) = l_k(2(x-x_j)/h_j - 1).
  \end{equation}
  We can then write
  \begin{equation}
    e_j(x) := e(x)|_{I_j} =  \sum_{k=0}^q \alpha_{j,k} l_{j,k}(x)  \text{ for } j=0, \dots, N-1.
  \end{equation}
  Since
  \begin{equation}
    \int_{I_j} l_{j,k}(x) l_{j,m}(x) \d x = 0 \Foreach k \neq m
  \end{equation}
  by the orthogonality condition on $e$ we can conclude $\alpha_{j,k} =
  0$ for $k=0, \dots, q-1$ and hence
  \begin{equation}
    e_j(x) = \alpha_{j,q} l_{j,q}(x).
  \end{equation}
  Now making use of the second condition
  \begin{equation}
    \begin{split}
      \avg{e}_j &= \frac 1 2 \qp{e(x_j^-) + e(x_j^+)}
      \\
      &= \frac 1 2 \qp{e_{j-1}(x_j^-) + e_j(x_j^+)}
      \\
      &= \frac 1 2 \qp{
        \alpha_{j-1,q} l_{j-1,q}(x_j) + \alpha_{j,q} l_{j,q}(x_j)}
      \\
      &= \frac 1 2 \qp{
        \alpha_{j-1,q} + (-1)^q \alpha_{j,q} },
    \end{split}
  \end{equation}
  through properties of the Legendre polynomials. Taking into account
  all edge contributions, this yields a linear system for the
  coefficients $\alpha_{j,q}$,
  \begin{equation}
    \begin{bmatrix}
      (-1)^q & 0 & 0 & \dots & 1
      \\
      1 & (-1)^q & 0 & \dots & 0
      \\
      0 & 1 & (-1)^q & \dots & 0
      \\
      \vdots & \vdots & \vdots & \ddots & \vdots
      \\
      0 & 0 & 0 & \dots & (-1)^q
    \end{bmatrix}
    \begin{bmatrix}
      \alpha_{0,q}\\
      \alpha_{1,q}\\
      \alpha_{2,q}\\
      \vdots\\
      \alpha_{N-1,q}
    \end{bmatrix}
    =
    \begin{bmatrix}
      v(x_0) - \cT_h(x_0^-)\\
      v(x_1) - \cT_h(x_1^-)\\
      v(x_2) - \cT_h(x_2^-)\\
      \vdots\\
      v(x_{N-1}) - \cT_h(x_{N-1}^-)\\
    \end{bmatrix}
    .
  \end{equation}
  This system is invertible only when $N$ is odd and $q$ is even
  \revise{which also ensures uniqueness of the projector $\cS_h(v)$}.

  Now solving this system, we have
  \begin{equation}
    \begin{bmatrix}
      \alpha_{0,q}\\
      \alpha_{1,q}\\
      \alpha_{2,q}\\
      \vdots\\
      \alpha_{N-1,q}
    \end{bmatrix}
    =
    \frac 1 2
    \begin{bmatrix}
      1 & 1 & -1 & \dots & 1 & -1
      \\
      -1 & 1 & 1 & \dots & -1 & 1
      \\
      1 & -1 & 1 & \dots & 1 & -1
      \\
      \vdots & \vdots & \vdots & \ddots & \ddots & \vdots
      \\
      1 & -1 & 1 & \dots & -1 & 1
    \end{bmatrix}
    \begin{bmatrix}
      v(x_0) - \cT_h(x_0^-)\\
      v(x_1) - \cT_h(x_1^-)\\
      v(x_2) - \cT_h(x_2^-)\\
      \vdots\\
      v(x_{N-1}) - \cT_h(x_{N-1}^-)\\
    \end{bmatrix}
    .
  \end{equation}
  
  From (\ref{eq:T-approx}) we know
  \begin{equation}
    \norm{ v(x_j) - \cT_h(x_j^-)}
    \leq
    C h^{q+1} \Norm{v}_{\sob{q+1}{\infty}(I_j)} \text{ for } j = 0,\dots, N-1.
  \end{equation}
  To conclude we invoke the results of \cite[Prop
    3.2]{BonaChenKarakashianXing:2013} that states
  \begin{equation}
    \norm{ v(x_i) - \cT_h(x_i^-) - v(x_{i+1}) + \cT_h(x_{i+1}^-)}
    \leq C h^{q+2} \Norm{v}_{\sob{q+2}{\infty}(I_i \cup I_{i+1})} \text{ for } i = 0,\dots, N-2.
  \end{equation}
  Appropriately extending the result over the periodic boundary we have
  \begin{equation}
    \norm{\alpha_{j,q}} \leq C h^{q+1} \norm{u}_{\sob{q+2}{\infty}(S^1)},
  \end{equation}
  by (\ref{eq:T-approx}). This yields $\leb{\infty}(S^1)$ control on
  the error and $\leb{2}(S^1)$ control as a consequence. For the
  gradient bound, note through the definition of $\cG_h$ we have
  \begin{equation}
    \begin{split}
      \Norm{\cP_h \qp{v_x} - \cG_h(\cS_h(v))}_{\leb{2}(S^1)}
      &=
      \sup_{\phi\in\leb{2}, \Norm{\phi}\leq 1}
      \int_{S^1} \qp{\cP_h \qp{v_x} - \cG_h(\cS_h(v))} \phi
      \\
      &=
      \sup_{\phi\in\leb{2}, \Norm{\phi}\leq 1}
      \int_{S^1} \qp{v_x - \cG_h(\cS_h(v))} \cP_h(\phi)
      \\
      &=
      \sup_{\phi\in\leb{2}, \Norm{\phi}\leq 1}
      -\int_{S^1} \qp{v - \cS_h(v)} \qp{\cP_h(\phi)}_x \d x
      +
      \sum_{j=0}^{N-1} \avg{v - \cS_h(v)}_j\jump{\cP_h(\phi)}_j
      \\
      &=0,
    \end{split}
  \end{equation}
  by the definition of $\cS_h$. Hence $\cP_h \qp{v_x} =
  \cG_h(\cS_h(v))$ and the result follows through the approximation
  properties of the $\leb{2}$ projection.
\end{proof}

\begin{Remark}[Restrictions of Lemma \ref{lem:S}]
  The a priori analysis that follows has restrictions that stem from
  defining the projection operator $\cS_h$, the mesh-size should be
  uniform, although it is possible to relax this condition, the
  polynomial degree should be even and the number of grid points
  odd. The reason is to have access to an operator that is consistent
  with our discretisation of the discrete derivative
  operator. Numerically, \revise{we observe optimal convergence of the
    method regardless of the polynomial degree and the number of grid
    points.}

  Note that similar conditions are required for the dispersive
  projection operator in \cite{BonaChenKarakashianXing:2013}. These
  conditions can be removed using an L-dG approach
  \cite{YanShu:2002}. We have chosen to pursue this approach to enable
  us to use elliptic reconstruction techniques for the a posteriori
  analysis following in Section \ref{sec:apos}.
\end{Remark}

\begin{Definition}[Mesh dependent norms]
  Let $V(h):= \fes_q + \sobh2(S^1)$. We define two mesh dependent
  $\sobh1$- like norms as
  \begin{equation}
    \begin{split}
      \denorm{W}^2
      &:=
      \sum_{j=0}^{N-1} \qp{\Norm{W_x}_{\leb{2}(I_j)}^2 + \Norm{h W_{xx}}_{\leb{2}(I_j)}^2 + \avg{h}_j^{-1} \jump W_j^2 }
      \\
      \enorm{W}^2
      &:=
      \sum_{j=0}^{N-1} \qp{\Norm{W_x}_{\leb{2}(I_j)}^2 + \avg{h}_j^{-1} \jump W_j^2 }.
    \end{split}
  \end{equation}
  When the penalty parameter, $\sigma$, is chosen large enough, the
  interior penalty bi-linear form (\ref{def:ip}) is coercive over
  $\fes_q$ and continuous over $V(h)$ with respect to the norm
  $\denorm{\cdot}$, that is
  \begin{equation}
    \begin{split}
      C_c \denorm{W}^2 &\leq \bih{W}{W} \Foreach W\in\fes_q
      \\
      \bih{w}{\phi} &\leq C_b \denorm{w}\denorm{\phi} \Foreach w,\phi \in V(h).
    \end{split}
  \end{equation}
  For $W\in\fes_q$ the two norms are equivalent and in particular
  \begin{equation}
    \enorm{W} \leq \denorm{W}.
  \end{equation}
  We will make use of $\denorm{\cdot}$ for a priori analysis and
  $\enorm{\cdot}$ and $\denorm{\cdot}$ for a posteriori analysis.
\end{Definition}

\begin{Lemma}[Inconsistent Ritz projector $\ritz$ and error control]
  \label{lem:R}
  Let the conditions in Lemma \ref{lem:S} hold. For $u\in
  \sob{q+3}{\infty}(S^1)$, let $\ritz(u)\in \fes_q$
  satisfy
  \begin{equation}
    \label{eq:utwidle}
    \bih{\ritz(u)}{\Phi} + \int_{S^1} \ritz(u)\Phi \d x
    =
    \bih{u}{\Phi}
    +
    \int_{S^1} \qp{u + v - \cS_h(v)} \Phi \d x
    .
  \end{equation}
  Then, for $h$ small enough
  \begin{equation}
    \label{eq:err-u}
    \Norm{u - \ritz(u)}_{\leb{2}(S^1)}
    +
    h \denorm{u - \ritz(u)}
    \leq
    C h^{q+1} {\Norm{u}_{\sob{q+3}{\infty}(S^1)}}.
  \end{equation}
\end{Lemma}
\begin{proof}
  To show (\ref{eq:err-u}) we note that through the definition
  (\ref{eq:utwidle}) we have the orthogonality result
  \begin{equation}
    \bih{\ritz(u) - u}{\Phi} + \int_{S^1} \qp{\ritz(u) - u}\Phi \d x
    =
    \int_{S^1} \qp{v - \cS_h (v)} \Phi \d x
    \Foreach \Phi \in \fes_q.
  \end{equation}
  Hence we have, for any $W\in\fes_q$
  \begin{equation}
    \begin{split}
      C_c\denorm{W - \ritz(u)}^2 + \Norm{W - \ritz(u)}_{\leb{2}(S^1)}^2
      &\leq      
      \bih{W - \ritz(u)}{W - \ritz(u)} + \int_{S^1} \qp{W - \ritz(u)}\qp{W - \ritz(u)} \d x
      \\
      &=
      \bih{W - u}{W - \ritz(u)} + \int_{S^1} \qp{W - u}\qp{W - \ritz(u)} \d x
      \\
      &\qquad +
      \bih{u - \ritz(u)}{W - \ritz(u)} + \int_{S^1} \qp{u - \ritz(u)}\qp{W - \ritz(u)} \d x
      \\
      &=
      \bih{W - u}{W - \ritz(u)} + \int_{S^1} \qp{W - u}\qp{W - \ritz(u)} \d x
      \\
      &\qquad +
      \int_{S^1} \qp{\cS_h(v) - v} \qp{W - \ritz(u)} \d x
      \\
      &\leq
      \frac 1 2 \bigg(
      \frac{C_b^2}{C_c}\denorm{W - u}^2
      \revise{
        +
        C_c\denorm{W - \ritz(u)}^2
        + \Norm{W - u}_{\leb{2}(S^1)}^2 }
      \\
      &\qquad + 2\Norm{W - \ritz(u)}_{\leb{2}(S^1)}^2
      +
      \Norm{\cS_h(v) - v}_{\leb{2}(S^1)}^2 
      \bigg).
    \end{split}
  \end{equation}
  Thus, choosing $W = \cP_h(u)$ and using approximation properties of
  the $\leb{2}$ projector as well as the bound from Lemma \ref{lem:S}
  we have
  \begin{equation}
    \denorm{W - \ritz(u)}^2
    \leq
    Ch^{2q}\qp{ \norm{u}_{\sobh{q+1}(S^1)}^2 + \Norm{v}_{\sob{q+1}{\infty}(S^1)}^2},
  \end{equation}
  and hence the $V(h)$ norm bound follows from the triangle inequality
  and the definition of $v$. To show the $\leb{2}$ bound, let
  $z\in\sobh 2(S^1)$ solve the dual problem
  \begin{equation}
    \begin{split}
      -z_{xx} + z &= u - \ritz(u),
    \end{split}
  \end{equation}
  then elliptic regularity guarantees that
  \begin{equation}
    \label{eq:ell-reg}
    \norm{z}_{\sobh2 (S^1)} \leq C \Norm{u - \ritz(u)}_{\leb{2}(S^1)}.
  \end{equation}
  Hence, for any $Z\in \fes_q$
  \begin{equation}
    \begin{split}
      \Norm{u - \ritz(u)}_{\leb{2}(S^1)}^2
      &=
      \int_{\W} \qp{u - \ritz(u)}^2 \d x
      \\
      &=
      \int_{\W}\qp{- z_{xx} + z}\qp{u - \ritz(u)} \d x
      \\
      &=
      \bih{z}{u - \ritz(u)}
      +
      \int_{S^1} z\qp{u-\ritz(u)} \d x
      \\
      &=
      \bih{z-Z}{u - \ritz(u)}
      +
      \int_{S^1} \qp{z-Z}\qp{u-\ritz(u)}
      +
      \qp{\cS_h(v) - v} Z \d x,
    \end{split}
  \end{equation}
  \revise{by the quasi-orthogonality result (\ref{eq:utwidle}).}
  Making use of the orthogonality of $\cS_h$, we choose
  $Z\in\fes_{q-1}$ as the $\leb{2}$ orthogonal projector of $z$ and
  find, by Cauchy-Schwarz 
  \begin{equation}
    \begin{split}
      \Norm{u - \ritz(u)}_{\leb{2}(S^1)}^2
      &\leq
      C_b \denorm{z-Z} \denorm{u - \ritz(u)}
      +
      \Norm{z-Z}_{\leb{2}(S^1)}\Norm{u-\ritz(u)}_{\leb{2}(S^1)}
      \\
      &\leq
      Ch\norm{z}_{2} \denorm{u - \ritz(u)}
      +
      Ch^2\norm{z}_2\Norm{u-\ritz(u)}_{\leb{2}(S^1)}
      \\
      &\leq
      Ch \Norm{u - \ritz(u)}_{\leb{2}(S^1)} \denorm{u - \ritz(u)}
      +
      Ch^2 \Norm{u-\ritz(u)}_{\leb{2}(S^1)}^2
    \end{split}
  \end{equation}
  using the elliptic regularity result (\ref{eq:ell-reg}) and
  approximation properties of the $\leb{2}$ projector. Hence
  \begin{equation}
    \qp{1-Ch^2}      \Norm{u - \ritz(u)}_{\leb{2}(S^1)}
    \leq
    Ch \denorm{u - \ritz(u)}
  \end{equation}
  as required for $h$ small enough.
\end{proof}

\begin{Theorem}[A priori bound - linear case]
  Suppose $f(u) = \frac 1 2 u^2$, in this case the PDE (\ref{thorder})
  is linear and given by
  \begin{equation}
    u_t - u_x + u_{xxx} = 0.
  \end{equation}
  Assume that the solution of (\ref{thorder}) $u \in
  \sob{q+4}{\infty}(S^1)$ and $u_t \in \sob{q+3}{\infty}(S^1)$ and let
  $U$ solve (\ref{sds}) for some $q$ even and $N$ odd. Then, for $t\in
      [0,T]$ and $h$ small enough,
  \begin{equation}
    \label{eqn:defocusapriori}
    \begin{split}
      \frac 1 2
      \denorm{\qp{u - U}(t)}^2 + \frac 1 2 \Norm{\qp{u - U}(t)}_{\leb{2}(S^1)}^2
      &\leq
      \exp{\qp{t}}\bigg(
      \frac 1 2\denorm{\qp{u - U}(0)}^2 + \frac 1 2\Norm{\qp{u - U}(0)}_{\leb{2}(S^1)}^2
      \\
      &\qquad +
      Ch^{2q}  \int_0^t
      \Norm{u_t}_{\sob{q+3}{\infty}(S^1)}^2
        +
        \Norm{u}_{\sob{q+4}{\infty}(S^1)}
        \d s
      \bigg).
    \end{split}
  \end{equation}
\end{Theorem}

\begin{proof}
  We begin by noting that, since $f'(u) = u$, in Lemma
  \ref{lem:drrer}, hence
  \begin{equation}
    \label{lem:pf1}
    \ddt F_{1,h}(\theta^u)
    =
    \int_{S^1}  \theta^u \Ru - \Rv \cG_h(\theta^v) \d x
    +
    \bih{\theta^u}{\Ru}
    .
  \end{equation}
  Observe that the term $\cG_h(\theta^v)$ is not controllable in
  $F_{1,h}(\theta^u)$ and also will not be of an optimal order. It is
  prudent for fixed $\qp{U, V}$ to choose $\qp{\widetilde U,
    \widetilde V}$ such that $\Rv = 0$. This then constrains choices
  for the pair $\qp{\widetilde U, \widetilde V}$. We pick $\widetilde
  V = \cS_h(v)$ and then choose $\widetilde U = \ritz(u)$. This choice
  ensures that the perturbed equations
  \begin{equation}
    \begin{split}
      \int_{S^1} \widetilde U_t \Phi
      + \cG_h(\widetilde V)\Phi \d x
      &=
      \int_{S^1} -\Ru \Phi \d x \quad \forall \ \Phi \in \rV_q \\
      \int_{S^1} \widetilde V \Psi
      + f'(\widetilde U) \Psi \d x
      + \bih{\widetilde U}{\Psi}
      &=\int_{S^1} -\Rv \Psi \d x \quad \forall \ \Psi \in \rV_q
      ,
    \end{split}
  \end{equation}
  are satisfied with
  \begin{equation}
    \begin{split}
      \Ru &= u_t - \widetilde U_t  + v_x - \cG_h(\widetilde V) 
      \\
      \Rv &= 0.
    \end{split}
  \end{equation}
  Substituting this into (\ref{lem:pf1}) we have 
    \begin{equation}
    \label{lem:pf2}
    \ddt F_{1,h}(\theta^u)
    =
    \int_{S^1} \theta^u \Ru \d x
    +
    \bih{\theta^u}{\Ru}
    .
    \end{equation}
    Now, through Cauchy's inequality we see 
    \begin{equation}
    \label{lem:pf3}
    \ddt F_{1,h}(\theta^u)
    \leq
    \frac 1 2
    \qp{
      \Norm{\theta^u}_{\leb{2}(S^1)}^2 + C_c \denorm{\theta^u}^2
      +
      \Norm{\Ru}_{\leb{2}(S^1)}^2 +
      \frac{C_b^2}{C_c}\denorm{\Ru}^2
    }
    .
    \end{equation}
    Now since
    \begin{equation}
      \frac{C_c}2 \denorm{\theta^u}^2 + \frac 1 2 \Norm{\theta^u}_{\leb{2}(S^1)}^2
      \leq
      F_{1,h}(\theta^u),
    \end{equation}
    Gronwall's inequality implies
    \begin{equation}
      \begin{split}
        C_c \denorm{\theta^u(t)}^2 + \Norm{\theta^u(t)}_{\leb{2}(S^1)}^2
        &\leq
        \exp{\qp{t}}
        \bigg(
          C_c \denorm{\theta^u(0)}^2 + \Norm{\theta^u(0)}_{\leb{2}(S^1)}^2
          \\
          &\qquad \qquad +
          \int_0^t
          \Norm{\Ru(s)}_{\leb{2}(S^1)}^2 + \frac{C_b^2}{C_c}\denorm{\Ru(s)}^2
          \d s
        \bigg).
      \end{split}
    \end{equation}    
    It remains to bound the term $\Ru$. We do this by splitting into
    two components and controlling them individually. First note that
    since we are in a semi discrete setting, Lemma \ref{lem:R} yields
    \begin{equation}
      \Norm{u_t - \ritz(u_t)}_{\leb{2}(S^1)}
      \leq
      Ch^{q+1} \Norm{u_t}_{\sob{q+3}{\infty}(S^1)}.
    \end{equation}
    Further, Lemma \ref{lem:S} immediately gives
    \begin{equation}
      \Norm{v_x - \cG_h(\cS_h(v))}_{\leb{2}(S^1)} \leq
      C h^{q+1} \Norm{u}_{\sob{q+4}{\infty}(S^1)},
    \end{equation}
    hence
    \begin{equation}
      \Norm{\Ru}_{\leb{2}(S^1)}^2 + C_b^2\denorm{\Ru}^2
      \leq
      C h^{2q} \qp{\Norm{u_t}_{\sob{q+3}{\infty}(S^1)}^2
        +
        \Norm{u}_{\sob{q+4}{\infty}(S^1)}^2
      },
    \end{equation}
    as required.
\end{proof}

\begin{Remark}[A priori bound - nonlinear case]
  In the nonlinear case, for $f(u) = \alpha u^m$, for $m>2$, our
  problem is given by
  \begin{equation}
    u_t - \alpha \qp{u^m}_x + u_{xxx} = 0.
  \end{equation}
  An optimal bound for this using the methodology proposed above
  requires the appropriate handling of discontinuous Galerkin
  approximations of the associated Emden-Fowler type elliptic problem
  \begin{equation}
    -u_{xx} + \abs{u}^{m-2} u = f,
  \end{equation}
  which is discussed in \cite{self:quasinorm}. It should be noted that
  optimal a priori control of approximations to this problem even in
  the energy norm are not trivial for $m> 2$.
\end{Remark}

\section{A posteriori analysis}
\label{sec:apos}

In this section, we give an a posteriori analysis of the semi discrete
scheme posed in Section \ref{sec:dis}. We proceed along similar lines
to the a priori analysis in that we examine solutions of perturbed
equations, taking account of different effects errors induced will
have. The difference being, in this section we make use of the
stability framework of the underlying PDE.

\begin{Lemma}
  \label{lem:pertapost}
  Let $u \in \cont 1([0,T], \sobh3(S^1))$ be a strong solution to
  (\ref{thorder}) and suppose $\widetilde u\in \cont 1([0,T],
  \sobh3(S^1))$ satisfies the problem
  \begin{equation}
    \begin{split}
      \widetilde u_t - f'(u)_x + f'(u - \widetilde u)_x + \widetilde u_{xxx} &= -\Rapost,
      \\
      \widetilde u(x,0) &= \widetilde u^0(x)
    \end{split}
  \end{equation}
  for some $\Rapost \in \leb 2(S^1)$. Then, with $\rho := u - \widetilde u$
  \begin{equation}
    \ddt \qp{\int_{S^1} \frac 1 2 {\rho_x^2} + f(\rho) \d x}
    =
    \int_{S^1}\qp{-\rho_{xx} + f'(\rho)} \Rapost \d x
  \end{equation}  
\end{Lemma}

\begin{proof}
  To begin, we note that $\rho = u - \widetilde u$ satisfies the
  error equation
  \begin{equation}
    \label{eq:rho-eq}
    \rho_t - f'(\rho)_x + \rho_{xxx} = \Rapost.
  \end{equation}
  Then, explicitly computing the time derivative we have
  \begin{equation}
    \begin{split}
      \ddt F_1(\rho)
      &=
      \int_{S^1} \rho_x \rho_{xt} + f'(\rho)\rho_t \d x
      \\
      &=
      \int_{S^1} -\rho_{xx} \rho_{t} + f'(\rho) \rho_t \d x
      .
    \end{split}
  \end{equation}
  Making use of (\ref{eq:rho-eq}) we see
  \begin{equation}
    \begin{split}
      \int_{S^1} -\rho_{xx} \rho_{t} + f'(\rho) \rho_t \d x
      &=
      \int_{S^1} \qp{-\rho_{xx} + f'(\rho)}
      \qp{\Rapost - \rho_{xxx} + f'(\rho)_x} \d x
      \\
      &=
      \revise{
      \int_{S^1} \qp{-\rho_{xx} + f'(\rho)}
      \qp{\qp{ - \rho_{xx} + f'(\rho)}_x + \Rapost \d x }
      }
      \\
      &=
      \int_{S^1} \qp{-\rho_{xx} + f'(\rho)}
      \Rapost \d x
      ,
    \end{split}
  \end{equation}  
  as required.
\end{proof}

\begin{Definition}[Orthogonal decomposition]
  \label{def:orthog}
  To make use of the stability framework of the PDE, we split the
  numerical solution, $U\in\fes_q$ into a continuous component and a
  discontinuous component
  \begin{equation}
    U = U^c + U^d,
  \end{equation}
  where $U^c \in\cont{1}(0,T; \fes_q\cap \cont{0}(S^1))$ and $U^d := U
  - U^c \in\cont{1}(0,T; \fes_q)$. We decompose in such a way that
  \begin{equation}
    \label{eq:orth-decom}
    \bih{U^d}{\Phi} = 0 \Foreach \Phi\in\fes_q\cap \cont{0}(S^1),
  \end{equation}
  as in \cite{dedner2019residual}. We will make use of the following
  result, the proof of which can be found in
  \cite{HoustonPerugiaSchotzau:2004}.
\end{Definition}

\begin{Proposition}[Bound on nonconforming term]
  Let $U=U^c + U^d$ be the orthogonal decomposition from Definition
  \ref{def:orthog} then we have
  \begin{equation}
    \denorm{U^d}^2 \leq C \sum_{j=0}^{N-1} \avg{h}_j^{-1} \jump{U}_j^2,
  \end{equation}
  where $C$ depends only upon $q$ and the $\max_j \frac{h_j}{h_{j+1}}$.
\end{Proposition}

\renewcommand{\D}{\cD_h}

\begin{Theorem}[Discrete reconstruction operator $\D$]
  \label{eq:discrete-reconstruction}
  For each $W\in\fes_q$ there exists a unique $\D(W)\in
  \fes_{q+1}\cap{\cont{0}(S^1)}$ such that for $j=0, \dots, N-1$
  \begin{equation}
    \begin{split}
      \int_{S^1} \D (W)_x \Phi \d x
      &=
      \int_{S^1} \cG_h(W) \Phi \d x \Foreach\Phi\in\fes_{q} , \qquad t>0
      \\
      \D(W)(x_j^+) &= \avg{W}_j \Foreach j = 0,\dots N-1.
    \end{split}
  \end{equation}
  Moreover, $\D$ satisfies
  \begin{gather}
    \Norm{W - \D(W)}_{\leb{2}(S^1)}^2
    \leq C \sum_{j=0}^{N-1} \avg{h}_j \jump{W}_j^2.
  \end{gather}
\end{Theorem}

\begin{proof}
  Fix $W\in \fes_q$, then a candidate $\Psi\in\fes_q$ given by
  \begin{equation}
    \int_{x_j}^{x_{j+1}} \Psi \Phi \d x = \int_{x_j}^{x_{j+1}}
    \cG_h(W)\Phi \d x
    \revise{\Foreach j = 0,\dots N-1}
  \end{equation}
  exists through Riesz Representation Theorem. Now, with
  $\D(W)\in\fes_{q+1}$ as the anti spatial derivative of $\Psi$, the
  constant of integration can be chosen such that the boundary
  condition is satisfied showing existence and uniqueness.
  
  To show $\D(W)\in\cont{0}(S^1)$ we fix $j$ and consider

  \begin{equation}
    \Phi(x)
    =
    \begin{cases}
      1 \text{ if } x\in ({x_j},{x_{j+1}})\\
      0 \text{ otherwise.}
    \end{cases}
  \end{equation}
  Now note that by definition, \revise{and the periodic boundary
    conditions, we have}
  \begin{equation}
    \begin{split}
      0 &= \int_{x_j}^{x_{j+1}} \D (W)_x \Phi 
      -
      \cG_h(W) \Phi \d x
      \\
      &=
      \D (W)(x_{j+1}^-) -  \D (W)(x_{j}^+)
      -
      \avg{W}_{j+1} + \avg{W}_{j}
      \\
      &=
      \D (W)(x_{j+1}^-)
      -
      \avg{W}_{j+1}.
    \end{split}
  \end{equation}
  This ensures $\D(W)$ is continuous over the edge and, since $j$ was
  arbitrary, globally. To show the approximation properties we note
  that $\D(W) - W \in\fes_{q+1}$ and, in particular,
  \begin{equation}
    \label{eq:D-orthog}
    \int_{x_j}^{x_{j+1}} \qp{\D(W) - W}\Phi \d x = 0 \Foreach \Phi\in\fes_{q-1}.
  \end{equation}
  We can then write the difference in terms of the Legendre
  polynomials defined in the Proof of Lemma \ref{lem:S}. In
  particular, due to the orthogonality condition (\ref{eq:D-orthog})
  we have
  \begin{equation}
    \qp{\D(W) - W}|_{I_j} = \sum_{k=0}^{q+1} \alpha_{j,k} l_{j,k}(x) = \alpha_{j,q} l_{j,q}(x) + \alpha_{j,q+1} l_{j,q+1}(x) \text{ for } j=0,\dots, N-1
  \end{equation}
  and endpoint conditions
  \begin{equation}
    \begin{split}
      \qp{\D(W) - W}(x_j^+) &= \avg{W}_j - W_j^+ = \frac 12 \jump{W}_j
      \\
      \qp{\D(W) - W}(x_{j+1}^-) &= \avg{W}_{j+1} - W_{j+1}^- = -\frac 12 \jump{W}_{j+1}.
    \end{split}
  \end{equation}
  Making use of the properties of the Legendre polynomials, we can write the linear system
  \begin{equation}
    \begin{split}
      \alpha_{j,q} (-1)^q + \alpha_{j,q+1} (-1)^{q+1} &= \frac 12 \jump{W}_j
      \\
      \alpha_{j,q} + \alpha_{j,q+1} &= -\frac 12 \jump{W}_j.
    \end{split}
  \end{equation}
  This can then be readily solved to show
  \begin{equation}
    \begin{split}
      \alpha_{j,q} &= \frac 14 \qp{(-1)^q \jump{W}_j - \jump{W}_{j+1}}
      \\
      \alpha_{j,q+1} &= \frac 14 \qp{(-1)^{q+1} \jump{W}_j - \jump{W}_{j+1}},
    \end{split}
  \end{equation}
  giving an explicit representation for $\qp{\D(W) - W}|_{I_j}$. Hence
  \begin{equation}
    \begin{split}
      \int_{x_j}^{x_{j+1}} \qp{\D(W) - W}^2 \d x
      &=
      \int_{x_j}^{x_{j+1}}  
      \alpha_{j,q}^2 l_{j,q}(x)^2 + \alpha_{j,q+1}^2 l_{j,q+1}(x)^2
      \d x
      \\
      &=  
      \alpha_{j,q}^2 \frac{h_j}{2q+1} + \alpha_{j,q+1}^2  \frac{h_j}{2q+3}
      \\
      &\leq
      \frac{h_j}{8\qp{2q+1}} \qp{\jump{W}_j^2 + \jump{W}_{j+1}^2}.
    \end{split}
  \end{equation}
  The result follows by summing over all elements with $C=
  \qp{4\qp{2q+1}}^{-1}$ if the mesh is uniform. In the nonuniform
  setting $C$ depends on $\max_j \frac{h_j}{h_{j+1}}$, the grading of
  the mesh.
\end{proof}

\begin{Definition}[Elliptic reconstruction]
  \label{def:recon}
  Let $\qp{U,V}\in\fes_q\times\fes_q$ be the semi-discrete approximation given by
  (\ref{sds}). Then, the elliptic reconstruction $\recon : \fes_q \to
  \sobh 3(S^1)$ is given by the solution of
  \begin{equation}
    \label{eq:recon}
    -\recon_{xx} + f'(\recon) = \D(V)
  \end{equation}
  with average value matching the discrete solution, that is
  \begin{equation}
    \int_{S^1} \recon - U \d x = 0.
  \end{equation}
\end{Definition}

\begin{Remark}[Inconsistent elliptic reconstruction]
  \label{rem:incon}
  The reconstruction $\recon$ is an inconsistent elliptic
  reconstruction of $U$ \cite{Makridakis:2003}. Indeed, $U$ is the
  finite element approximation of
  \begin{equation}
    -\recon_{xx} + f'(\recon) = V.
  \end{equation}
  The reason for defining $\recon$ as in Definition \ref{def:recon} is
  to ensure that the reconstruction has sufficient regularity to
  satisfy the perturbed PDE in Lemma \ref{lem:pertapost}.
\end{Remark}

\begin{Proposition}[Regularity bound for the reconstruction]
  \label{pro:regularity}
  The elliptic problem defining the reconstruction operator, $\recon$,
  in Definition \ref{def:recon} is well posed, moreover, thanks to
  elliptic regularity, we have
  \begin{equation}
    \Norm{\recon}_{\sobh{k+1}(S^1)}
    \leq
    C_{reg}
    \Norm{\D(V)}_{\sobh{k-1}(S^1)}
    \qquad \text{ for } k = 0, 1, 2.
  \end{equation}
\end{Proposition}

\begin{Lemma}[Reconstructed PDE]
  \label{lem:recon}
  The reconstruction given in Definition \ref{def:recon} satisfies
  \begin{equation}
    \label{eq:recon2}
    \recon_t - f'(u)_x + f'(u - \recon)_x + \recon_{xxx}
    = \Rapost,
  \end{equation}
  with
  \begin{equation}
    \Rapost = \qp{\recon - U}_t + f'(\recon)_x - f'(u)_x + f'(u - \recon)_x.
  \end{equation}
\end{Lemma}
\begin{proof}
  Since $\recon$ satisfies (\ref{eq:recon}) and the problem data
  $\D(V) \in \sobh1(S^1)$ it is clear that $\recon \in \sobh3(S^1)$
  and satisfies
  \begin{equation}
    \label{eq:recon1}
    -\recon_{xxx} + f'(\recon)_x = \D(V)_x = U_t,
  \end{equation} 
  upon writing the definition of $\D(V)$ point-wise. Now, substituting
  (\ref{eq:recon1}) into (\ref{eq:recon2}), we see
  \begin{equation}
    \begin{split}
      \recon_t - f'(u)_x + f'(u - \recon)_x + \recon_{xxx}
      &=
      \qp{\recon - U}_t + f'(\recon)_x - f'(u)_x + f'(u - \recon)_x,
    \end{split}
  \end{equation}
  as required.
\end{proof}

\begin{Hyp}[A posteriori control for the elliptic problem]
  \label{hyp:elliptic-aposteriori}
  We make the assumption that there exists an optimal order elliptic a
  posteriori estimate controlling the energy norm error. That is,
  there exists a functional $\eta$ depending only upon $U$ and the
  problem data such that
  \begin{equation}
    \begin{split}
      \qp{
        \Norm{U -\recon}_{\leb{m}(S^1)}^m
        +
        \enorm{U -\recon}^2
      }^{1/2}
      &\leq
      \eta(U, g, \sobh1(S^1)) 
      .
    \end{split}
  \end{equation}
\end{Hyp}

\begin{example}
  For $f(u) = \frac 1m u^m$ with $g := - \D(V)$ in
  \cite{self:quasinorm} it was shown that
  \begin{equation}
    \begin{split}
      \Norm{U -\recon}_{\leb{m}(S^1)}^m
      +
      \enorm{U -\recon}^2
      &\leq
      C\sum_{j=0}^{N-1}
      \bigg[
        \Norm{h \qp{g + U_{xx} - f'(U)}}_{\leb{2}(x_j, x_{j+1})}^2
        \\
        &\qquad +
        \avg{h}_j \jump{U_x}_j^2
        +
        \sigma \avg{h}_j^{-1}  \jump{U}_j^2
        +
        \avg{h}_j  \jump{V}_j^2
      \bigg]
    \end{split}
  \end{equation}
  satisfies Hypothesis \ref{hyp:elliptic-aposteriori}.

  Note that the inconsistency described in Remark \ref{rem:incon}
  induced by modifying the elliptic reconstruction is accounted for by
  the last term in this estimator.
\end{example}

\begin{Remark}[Alternative estimators]
  One of the strengths of the elliptic reconstruction methodology is
  the ability to use estimators that are not residual based. Indeed,
  recovery based a posteriori estimators have been widely used since
  their introduction by the engineering community in the 1980s. Their
  success is due to their simplicity of implementation, mild
  dependence of problem data and super-convergence properties. Work
  carried out on recovery estimators has reached a state of maturity
  for elliptic problems, see
  \cite{Ainsworth:2000,Bank:2003:2,Zienkiewicz:1987,LakkisPryer:2011a}
  and subsequent references. These estimators could also be used in the
  subsequent analysis.
\end{Remark}

\begin{Theorem}[A posteriori bound - linear case]
  \label{the:apost-linear}
  Suppose $f(u) = \frac 1 2 u^2$. Further, let $U$ solve (\ref{sds})
  and $\recon$ be the elliptic reconstruction from Definition
  \ref{def:recon}. Let the conditions of Lemma \ref{lem:pertapost}
  hold. Then, with $\rho:=u-\recon $, for $t\in [0,T]$,
  \begin{equation}
    \label{eq:aposteriori-linear}
    \begin{split}
      F_1(\rho(t))
      &\leq    
      \exp(t)
      \bigg(
      F_1(\rho(0))
      +
      \int_0^t
      \Norm{\qp{\recon_t - U_t}(s)}_{\leb{2}(S^1)}^2
      \\
      &\qquad \qquad \qquad \qquad +
      2\Norm{\qp{\recon_t - U_t^c}_x(s)}_{\leb{2}(S^1)}^2
      +
      2{C_a ^2 C_b^2} \denorm{U_t^d(s)}^2
      \d s
      \bigg).
    \end{split}
  \end{equation}
\end{Theorem}
\begin{proof}
  Since $f'(u) = u$, in Lemma \ref{lem:pertapost} $\Rapost = \recon_t - U_t$,
  hence
  \begin{equation}
    \label{eq:errorcont}
    \begin{split}
      \d_t F_1(\rho)
      &=
      \int_{S^1} \qp{-\rho_{xx} + \rho} \Rapost \d x
      \\
      &=
      \bih{\rho}{\recon_t - U_t}
      +
      \rho \qp{\recon_t - U_t} \d x.
    \end{split}
  \end{equation}
  Now making use of the orthogonal decomposition of $U = U^c + U^d$
  given in (\ref{eq:orth-decom}) we have
  \begin{equation}
    \begin{split}
      \bih{\rho}{\recon_t - U_t}
      &=
      \int_{S^1} \rho_x \qp{\recon_t - U_t^c}_x \d x
      -
      \bih{\rho}{U_t^d}
      \\
      &=
      \int_{S^1} \rho_x \qp{\recon_t - U_t^c}_x \d x
      -
      \bih{\rho - P}{U_t^d} \Foreach P \in\fes_q\cap \sobh1(S^1),
    \end{split}
  \end{equation}
  \revise{since $\jump{P}=0$.}  Choosing $P = I_h(\rho)$ as the
  Scott-Zhang interpolant of $\rho$, we have
  \begin{equation}
    \label{eq:clevertrick}
    \bih{\rho}{\recon_t - U_t}
    \leq
    \frac 1 4
    \Norm{\rho_x}_{\leb{2}(S^1)}^2
    +
    \Norm{\qp{\recon_t - U_t^c}_x}_{\leb{2}(S^1)}^2
    +
    {\epsilon} \denorm{\rho - I_h(\rho)}^2
    +
    \frac{C_b^2}{4\epsilon}\denorm{U_t^d}^2
    .
  \end{equation}
  Using the stability of the Scott-Zhang interpolant in mesh dependent
  norms, see \cite[Lemma 3.49 c.f.]{ErnGuermond:2004} we have
  \begin{equation}
    \denorm{\rho - I_h(\rho)} \leq C_a \Norm{\rho_x}_{\leb{2}(S^1)}.
  \end{equation}
  Hence, for any
  $\epsilon>0$
  \begin{equation}
    \label{eq:stabbound}
    \begin{split}
      \bih{\rho}{\recon_t - U_t}
      &\leq
      \frac 1 4
      \Norm{\rho_x}_{\leb{2}(S^1)}^2
      +
      \Norm{\qp{\recon_t - U_t^c}_x}_{\leb{2}(S^1)}^2
      +
      C_a^2 \epsilon \Norm{\rho_x}^2_{\leb{2}(S^1)}
      +
      \frac{C_b^2}{4 \epsilon} \denorm{U_t^d}^2
      \\
      &\leq
      \frac 1 2
      \Norm{\rho_x}_{\leb{2}(S^1)}^2
      +
      \Norm{\qp{\recon_t - U_t^c}_x}_{\leb{2}(S^1)}^2
      +
      {C_a ^2 C_b^2} \denorm{U_t^d}^2
    \end{split}
  \end{equation}
  with $\epsilon = 1/(4C_a^2)$. Substituting (\ref{eq:stabbound}) into
  (\ref{eq:errorcont}) we have 
  \begin{equation}
    \begin{split}
      \ddt F_1(\rho)
      \leq
      \frac 1 2
      \Norm{\rho_x}_{\leb{2}(S^1)}^2
      +
      \Norm{\qp{\recon_t - U_t^c}_x}_{\leb{2}(S^1)}^2
      +
      \frac 12 \Norm{\rho}_{\leb{2}(S^1)}^2
      +
      \frac 12 \Norm{\recon_t - U_t}_{\leb{2}(S^1)}^2
      +
      {C_a ^2 C_b^2} \denorm{U_t^d}^2
      .
    \end{split}
  \end{equation}
  The result follows from Gronwall's inequality.
\end{proof}

\begin{Corollary}[Computable a posteriori bound - linear
  case] \label{cor:apost-linear} 
  Let the conditions of Theorem \ref{the:apost-linear} hold. Then,
  with $e:=u-U$, for $t\in [0,T]$,
  \begin{equation}
    \label{eq:aposteriori-linear2}
    \begin{split}
      \frac12 \enorm{e(t)}^2 + \frac12 \Norm{e(t)}_{\leb{2}(S^1)}^2
      &\leq
      \eta(U(t), g, \sobh{1}(S^1))^2 \\
      & \qquad +
      2 
      \exp(t)
      \bigg(
      F_1(\rho(0))
      +
      \int_0^t
      2\eta(U_s(s), g_s, \sobh{1}(S^1))^2 \\
      & \qquad \qquad \qquad \qquad \qquad \qquad
      + \qp{2{C_a ^2 C_b^2}}\denorm{U_s^d(s)}^2
      \d s
      \bigg)
      \\
      &=:
      E_{st}(U(t))
      .
    \end{split}     
  \end{equation}
  \revise{Note that $E_{st}$ is the computed quantity in the numerical
    results in Section \ref{sec:numerics}.}
\end{Corollary}

\begin{proof}

  We begin by noting
  \begin{equation}
    \begin{split}
      \frac 14 \enorm{e}^2 + \frac14 \Norm{e}_{\leb{2}(S^1)}^2
      & \le
      F_1(\rho)
      + \frac12 \enorm{\theta}^2
      + \frac12 \Norm{\theta}_{\leb{2}(S^1)}^2
      ,
    \end{split}
  \end{equation}
  where $\rho:=u-\recon$ and $\theta:=\recon-U$. 
  Invoking Hypothesis \ref{hyp:elliptic-aposteriori} and applying
  Theorem \ref{the:apost-linear} we may conclude.

\end{proof}

\begin{Remark}[The linear vs nonlinear case]
  Notice that when $f(u) = \frac 1 2 u^2$ we have that $\widetilde u$
  solves a compatible KdV-like problem
  \begin{equation}
    \widetilde u_t - \widetilde u_x + \widetilde u_{xxx} = -\Rapost.
  \end{equation}
  
  This makes the analysis considerably simpler than for general
  $f$. For expositions sake in the rest of this section we will only
  consider the case $f(u) = \frac 1 4 u^4$, corresponding to the
  defocusing mKdV equation. The arguments for more general
  non-linearity are lengthy and we wish to highlight that our analysis
  allows for explicit control on the constants appearing in the a
  posteriori upper bound, something very challenging for nonlinear
  evolution problems.
\end{Remark}
\begin{Lemma}[A priori solution control]
  \label{lem:aprioriapost}
  Let $u$ solve (\ref{thorder}) and $\recon$ be the elliptic
  reconstruction given in Definition \ref{def:recon}. Then, with $f(u)
  = \frac 14 u^4$, and for an initial condition $u^0 \in
  \sobh{2}(S^1)$, we have the following
  \begin{gather}
    \label{eq:ap1}
    \Norm{\qp{u\recon}_{x}}_{\leb{\infty}(S^1)}
    \leq
    \Norm{\qp{u\recon}_{xx}}_{\leb{2}(S^1)} \leq C_4,
  \end{gather}
  with
  \begin{equation}
    \begin{split}
      C_4
      &:=
      C_{reg}\bigg(
      \Norm{\D(V)}_{\sobh{-1}(S^1)}
      \qp{
        \Norm{u^0_{xx}}_{\leb{2}(S^1)}^2
        +
        5 \Norm{u^0u^0_{x}}_{\leb{2}(S^1)}^2
        +
        \frac 1 2 \Norm{u^0}_{\leb{6}(S^1)}^6
      }^{1/2}
      \\
      &\qquad +
      \Norm{\D(V)}_{\leb{2}(S^1)}
      \qp{
        \Norm{u^0_{x}}_{\leb{2}(S^1)}^2
        +
        \Norm{u^0}_{\leb{4}(S^1)}^4
      }^{1/2}\bigg)
      .
    \end{split}
  \end{equation}
\end{Lemma}
\begin{proof}
  We
  begin by noting that with (\ref{thorder}) admits a further
  invariant. That is,
  \begin{equation}
    \ddt \int_{S^1} 2 u_{xx}^2 + 10 u^2 u_x^2 + u^6 \d x = 0.
  \end{equation}
  Hence we have
  \begin{equation}
    \label{eq:reg1}
    \Norm{u_{xx}}_{\leb{2}(\W)}
    \leq
    \qp{\Norm{u^0_{xx}}^2_{\leb{2}(\W)}
    +
    5 \Norm{u^0 u^0_x}^2_{\leb{2}(\W)}
    +
    \frac 1 2 \Norm{u^0}^6_{\leb{6}(S^1)}}^{1/2}.
  \end{equation}
  In addition, since
  \begin{equation}
    \ddt \int_{S^1} \frac 12 u_{x}^2 + \frac 14 u^4 \d x = 0,
  \end{equation}
  we have
  \begin{equation}
    \label{eq:reg2}
    \Norm{u_{x}}_{\leb{2}(\W)}
    \leq
    \qp{\Norm{u^0_{x}}_{\leb{2}(\W)}^2
    +   
    \frac 1 2 \Norm{u^0}_{\leb{4}(\W)}^4}^{1/2}.
  \end{equation}
  The result (\ref{eq:ap1}) follows from (\ref{eq:reg1}),
  (\ref{eq:reg2}) and Proposition \ref{pro:regularity} and an
  interpolation argument.
\end{proof}

\begin{Theorem}[A posteriori bound - nonlinear case $m=4$]
  \label{the:apost-nonlinear}
  Suppose $f(u) = \frac 1 4 u^4$. Further, let $U$ solve
  (\ref{sds}) and the conditions of Lemmas \ref{lem:pertapost} and
  \ref{lem:recon} hold.  Then, with $\rho:= u-\recon$, for $t\in
      [0,T]$,
  \begin{equation}
    \label{eq:aposteriori-nonlinear}
    \begin{split}
      F_1(\rho(t))
      &\leq    
      \exp(\cC_{4} t)
      \bigg(
      F_1(\rho(0))
      +
      \int_0^t
      \frac 1 4
      \Norm{\qp{\recon_t-U_t}(s)}_{\leb{m}(S^1)}^4
      \\
      &\qquad \qquad\qquad \qquad\qquad \qquad
      +
      \Norm{\qp{\recon_t - U_t^c}_x(s)}_{\leb{2}(S^1)}^2
      +
      C_a^2C_b^2
      \denorm{U_t^d(s)}^2
      \d s
      \bigg),
    \end{split}
  \end{equation}
  with
  \begin{equation}
    \begin{split}
      \cC_{4} &=
      \max
      \qp{\frac 1 2 + \frac {15} 2 C_4
        ,
        \frac 3 4 + \frac 9 4 C_4
      }
      .
    \end{split}
  \end{equation}
\end{Theorem}
\begin{proof}
  For $m=4$, since $f'(u) = u^3$ we have
  \begin{equation}
    \begin{split}
      \Rapost
      &=
      \recon_t - U_t + f'(\recon)_x - f'(u)_x + f'(u - \recon)_x
      \\
      &= 
      \recon_t - U_t - \qp{3u\recon \qp{u - \recon}}_x.
    \end{split}
  \end{equation}
  Initial inspection of the form of $\Rapost$ indicates the resultant
  bound should not be optimal. Surprisingly, this is not the case as
  the extra derivative can be ``hidden'' by requiring regularity of $u$
  and $\recon$, which we have already quantified in a
  computational fashion.

  Through Lemma \ref{lem:pertapost}
  \begin{equation}
    \label{eq:nlapos1}
    \begin{split}
      \d_t F_1(\rho)
      &=
      \int_{S^1} \qp{-\rho_{xx} + \rho^3} \Rapost \d x
      \\
      &=
      \underbrace{\bih{\rho}{\recon_t - U_t} }_{\cI_1}
      -
      \underbrace{3\bi{\rho}{\qp{u \recon \rho}_x}}_{\cI_2}
      +
      \int_{S^1}
      \underbrace{\rho^3 \qp{\recon_t - U_t}}_{\cI_3}
      -
      \underbrace{3\rho^3 \qp{u \recon \rho}_x}_{\cI_4} \d x.
    \end{split}
  \end{equation}
  We proceed to control these terms individually. To begin, arguing as
  in (\ref{eq:clevertrick}), making use of the orthogonal
  decomposition (\ref{eq:orth-decom}) with $P = I_h(\rho)$ as the
  Scott-Zhang interpolant of $\rho$
  \begin{equation}
    \label{eq:nlapos2}
    \begin{split}
      \cI_1 &= \bih{\rho}{\recon_t - U_t}
      \\
      &=
      \int_{S^1} \rho_x \qp{\recon_t - U_t^c}_x \d x - \bih{\rho - P}{U_t^d}
      \\
      &\leq
      \frac 1 2 \Norm{\rho_x}_{\leb{2}(S^1)}^2
      +
      \Norm{\qp{\recon_t - U_t^c}_x}^2_{\leb{2}(S^1)}
      +
      C_a^2 C_b^2 \denorm{U_t^d}^2.
    \end{split}
  \end{equation}
  Now, through expanding derivatives
  \begin{equation}
    \begin{split}
      3\qp{u\recon \rho}_{xx}
      =
      3\qp{u\recon}_{xx} \rho
      +
      6\qp{u\recon}_x \rho_x
      +
      3u\recon \rho_{xx},
    \end{split}
  \end{equation}
  and integrating by parts we see
  \begin{equation}
    \label{eq:pf321}
    \begin{split}
      \cI_2
      &=
      \int_{S^1}
      3\rho_x \qp{u \recon \rho}_{xx}
      \d x
      \\
      &=
      \int_{S^1}
      3\qp{u\recon}_{xx} \rho \rho_x
      +
      6\qp{u\recon}_x \rho_x^2
      +
      3u\recon \rho_{xx}\rho_x
      \d x.      
    \end{split}
  \end{equation}
  Now note that
  \begin{equation}
    \begin{split}
      \int_{S^1} 3u\recon \rho_{xx}\rho_x \d x
      &=
      - \int_{S^1} 3\qp{u\recon\rho_x}_x \rho_x \d x
      \\
      &=
      - 3\int_{S^1} \qp{u\recon}_x \rho_x^2 + u\recon\rho_{xx} \rho_x \d x.
    \end{split}
  \end{equation}
  Hence
  \begin{equation}
    \label{eq:pf123}
    \int_{S^1} 3u\recon \rho_{xx}\rho_x \d x
    =
    - \frac 32\int_{S^1} \qp{u\recon}_x \rho_x^2.
  \end{equation}
  Substituting (\ref{eq:pf123}) into (\ref{eq:pf321})
  \begin{equation}
    \cI_2
    =
    \int_{\W} 3 \qp{u \recon}_{xx} \rho \rho_x + \frac 9 2 \qp{u\recon}_x \rho_x^2
    \d x.
  \end{equation}
  Making use of H\"older's inequality \revise{and a Sobolev embedding}
  \begin{equation}
    \label{eq:nlapos3}
    \begin{split}
      \cI_2
      &\leq
      3 \Norm{\qp{u \recon}_{xx}}_{\leb{2}(S^1)}
      \Norm{\rho_x}_{\leb{2}(S^1)}
      \Norm{\rho}_{\leb{\infty}(S^1)}
      +
      \frac 9 2 \Norm{\qp{u \recon}_{x}}_{\leb{\infty}(S^1)}
      \Norm{\rho_x}_{\leb{2}(S^1)}^2
      \\
      &\leq
      3 C_4
      \Norm{\rho_x}_{\leb{2}(S^1)}^2 
      +
      \frac 9 2
      C_4
      \Norm{\rho_x}_{\leb{2}(S^1)}^2,
      \\
      &\leq
      \frac {15} 2
      C_4
      \Norm{\rho_x}_{\leb{2}(S^1)}^2
      ,     
    \end{split}
  \end{equation}
  by Lemma \ref{lem:aprioriapost}. The third term
  \begin{equation}
    \label{eq:nlapos4}
    \begin{split}
      \cI_3 &= \int_{S^1} \rho^3 \qp{\recon_t - U_t} \d x
      \\
      &\leq
      \Norm{\rho^3}_{\leb{\frac 43}(\W)} \Norm{\recon_t - U_t}_{\leb{4}(\W)}
      \\
      &\leq
      \Norm{\rho}_{\leb{4}(\W)}^3 \Norm{\recon_t - U_t}_{\leb{4}(\W)}
      \\
      &\leq
      \frac 3 4 \Norm{\rho}_{\leb{4}(\W)}^4
      +
      \frac 1 4 \Norm{\recon_t - U_t}_{\leb{4}(\W)}^4.
    \end{split}
  \end{equation} 
  The final term can be controlled by noticing
  \begin{equation}
    \begin{split}
      \cI_4
      =
      3\int_{S^1} \rho^3 \qp{u \recon \rho}_x \d x
      =
      3\int_{S^1} \frac 1 4\qp{\rho^4}_x u \recon + \rho^4 \qp{u \recon}_x \d x
      =
      \frac 9 4 \int_{S^1} \rho^4 \qp{u\recon}_x \d x.
    \end{split}
  \end{equation}
  Hence
  \begin{equation}
    \label{eq:nlapos5}
    \begin{split}
      \cI_4
      &\leq
      \frac 9 4
      \Norm{\qp{u\recon}_x}_{\leb{\infty}(S^1)}
      \Norm{\rho}_{\leb{4}(S^1)}^4
      \\
      &\leq
      \frac 9 4 C_4
      \Norm{\rho}_{\leb{4}(S^1)}^4,
    \end{split}
  \end{equation}
  by Lemma \ref{lem:aprioriapost}. Collecting the results
  (\ref{eq:nlapos2}), (\ref{eq:nlapos3}), (\ref{eq:nlapos4}) and
  (\ref{eq:nlapos5}), substituting into (\ref{eq:nlapos1}) we have
  \begin{equation}
    \begin{split}
      \d_t F_1(\rho)
      \leq     
      \cC_{4} F_1(\rho)
      +
      \Norm{\qp{\recon_t - U_t^c}_x}_{\leb{2}(S^1)}
      +
      \frac 1 4
      \Norm{\recon_t - U_t}_{\leb{4}(S^1)}^4
      +
      C_a^2 C_b^2 \denorm{U_t^d}^2.
    \end{split}
  \end{equation}
  The result follows from Gronwall's inequality.
\end{proof}

\begin{Corollary}[Computable a posteriori bound - nonlinear case, $m=4$]
  Let the conditions of Theorem \ref{the:apost-nonlinear} hold.  Then,
  with $e:=u-U$, for $t\in [0,T]$, and $m=4$ we have
  \begin{equation}
    \label{eq:aposteriori-nonlinear-full}
    \begin{split}
      \frac12 \enorm{e(t)}^2
      + \frac14 \Norm{e(t)}_{\leb{4}(S^1)}^4
      &\leq
      2\eta(U(t), g(t), \sobh{1}(S^1))^2
      \\
      &\qquad \qquad +
      8
      \exp(\cC_4 t)
      \bigg(
      F_1(\rho(0))
      +
      \int_0^t
      2\eta(U_t(s), g_s(s), \sobh{1}(S^1))^2 \\
      & \qquad \qquad \qquad \qquad \qquad \qquad \qquad \qquad 
      +
      \qp{2C_a^2C_b^2}\denorm{U_s^d(s)}^2
      \d s
      \bigg)
            \\
      &=:
      E_{st}(U(t))
      .
    \end{split}     
  \end{equation}
  \revise{Note that, by an abuse of notation, $E_{st}$ is the computed
    quantity in the numerical results in Section \ref{sec:numerics}.}
\end{Corollary}
\begin{proof}
  This proof follows an equivalent argument to that made in the proof
  of Corollary \ref{cor:apost-linear}.
\end{proof}

\begin{Remark}The bound for the point-wise in time $\leb{m}$ error, $m=2,4$
  appearing on the left-hand side of (\ref{eq:aposteriori-linear}) and
  (\ref{eq:aposteriori-nonlinear-full}), is tight only for very short
  times. As we will observe in Section \ref{sec:temp-and-numerics} on
  a uniform mesh of size $h\to 0$ the gradient term $\enorm{u-U} =
  \Oh(h^q)$, while $\Norm{\qp{u-U}(t)}_{\leb{m}(S^1)} = \Oh(h^{q+1})$.
\end{Remark}

\section{Temporal discretisation and numerical benchmarking}
\label{sec:temp-and-numerics}
\renewcommand{\dt}[1]{\tau_{#1}}
\renewcommand{\dx}[1]{h_{#1}}

Practically, a fully discrete approximation scheme is required for
implementation. For the readers convenience we will present an
argument for designing a fully discrete scheme. We consider a time
interval $[0,T]$ subdivided into a partition of $N$ consecutive
adjacent sub-intervals whose endpoints are denoted
$t_0=0<t_1<\ldots<t_{N}=T$.  The $n$-th time-step is defined as
${\tau_n := t_{n+1} - t_{n}}$.  We will consistently use the shorthand
$y^n(\cdot):=y(\cdot,t_n)$ for a generic time function $y$. We also
denote $\nplush{y} := \frac{1}{2}\qp{y^n + y^{n+1}}$.

We consider the temporal discretisation of (\ref{sds}) as follows:
Given $U^0$, for $n\in [0, N-1]$ find $U^{n+1}$ such that
\begin{equation}
  \label{eq:dis-mixed-system-fd}
  \begin{split}
    0 &= \int_{\rS^1} \qp{\frac{U^{n+1} - U^n}{\tau_n}  + \cG_h(V^{n+1})} \Phi \d x
    \\
    0 &= \int_{\rS^1} \qp{V^{n+1} + \frac{f(U^{n+1}) - f(U^n)}{U^{n+1} - U^n}} \Psi \d x
    + \bih{U^{n+1/2}}{\Psi}  \Foreach \qp{\Phi, \Psi} \in \fes_q
    \\
    U^0 &= \Pi^0 u^0
  \end{split}
\end{equation}
where $\Pi^0$ denotes the $\leb{2}$ orthogonal projector into
$\fes_q$.

\begin{Theorem}[Conservativity of the fully discrete scheme]
  \label{the:cons-fullydiscrete}
  Let $\{ U^n \}_{n=0}^N$ be the fully discrete scheme generated
  by (\ref{eq:dis-mixed-system-fd}), then we have that
  \begin{equation}
    \label{eq:fd-mass}
    \int_{S^1} U^n \d x = \int_{S^1} U^0 \d x \Foreach n \in [0, N]
  \end{equation}
  and
  \begin{equation}
    \label{eq:fd-energy}
    \frac 1 2\bih{U^n}{U^n} + \int_{S^1} f(U^n) \d x
    =
    \frac 1 2\bih{U^0}{U^0} + \int_{S^1} f(U^0) \d x \Foreach n \in [0, N].
  \end{equation}
\end{Theorem}
\begin{proof}
  To show (\ref{eq:fd-mass}) it suffices to show that
  \begin{equation}
    \int_{S^1} U^{n+1} \d x = \int_{S^1} U^{n} \d x
  \end{equation}
  and then the result follows inductively. To this end, choosing $\Phi
  = 1$ in (\ref{eq:dis-mixed-system-fd})
  \begin{equation}
    \begin{split}
      0
      &=
      \int_{\rS^1} \qp{\frac{U^{n+1} - U^n}{\tau_n}  + \cG_h(V^{n+1})} \d x
      \\
      &=
      \int_{\rS^1} \frac{U^{n+1} - U^n}{\tau_n} \d x,
    \end{split}
  \end{equation}
  by the definition of $\cG_h$. To see (\ref{eq:fd-energy})
  \begin{equation}
    \begin{split}
      \frac 12\bih{U^{n+1}}{U^{n+1}} -& \frac 12\bih{U^n}{U^n}
      + \int_{S^1} f(U^{n+1}) - f(U^n) \d x
      \\
      &=
      \bih{U^{n+1} - U^n}{\nplush{U}}
      +
      \int_{S^1} \frac{f(U^{n+1}) - f(U^n)}{U^{n+1} - U^n} \qp{U^{n+1} - U^n} \d x
      \\
      &=
      -\int_{S^1} V^{n+1} \qp{U^{n+1} - U^n} \d x
      \\
      &=
      \int_{S^1} \tau V^{n+1} \cG_h(V^{n+1}) \d x = 0,
    \end{split}
  \end{equation}
  using the second equation of (\ref{eq:dis-mixed-system-fd}) and the
  skew-symmetry of $\cG_h$, concluding the proof.
\end{proof}

\begin{Remark}[Structure of the temporal
  discretisation] \label{rem:tstruct}
  
  The temporal discretisation given in (\ref{eq:dis-mixed-system-fd})
  is \emph{not} a Runge-Kutta method unless the problem is linear. It
  resembles a midpoint discretisation and is formally of second order,
  however the treatment of the non-linearity is different. Although
  construction of higher order methods is possible they become very
  complicated to write down so we will not press this point
  here. Further to the method of lines dG-difference scheme we
  propose, other discretisation methods are indeed possible.  The
  spatial discretisation can be coupled to space-time Galerkin
  procedures, using, for example a continuous Galerkin method in time
  to guarantee conservativity. It is even possible to make use of
  hybrid dG-cG approaches to construct flexible adaptive space-time
  schemes making use of recovered elements
  \cite{GeorgoulisPryer:2018}.
  
\end{Remark}
\begin{Remark}[Conservation of other invariants] \label{rem:tsym}
  \label{rem:cons}
  This discretisation does not lend itself to conservation of other
  invariants, for example even the quadratic invariant $F_0$ is not
  conserved under this scheme. A class of Runge-Kutta methods that are
  able to exactly conserve all quadratic invariants are the
  Gauss-Radau family, this is because they are \emph{symplectic}. When
  one considers higher order invariants, it seems that schemes must be
  designed individually and there seems to be no class that can
  exactly conserve all.
\end{Remark}

\section{Numerical experiments} \label{sec:numerics}

In this section we illustrate the performance of the method proposed
through a series of numerical experiments. The brunt of the
computational work has been carried out using Firedrake
\cite{Firedrake:2017}. We employ a Gauss quadrature of order $4q$,
where $q$ is the degree of the finite element space, to minimise
quadrature error introduced into the implementation. Indeed, at this
degree all integrals are performed exactly with the exception of the
projection of the initial condition. The nonlinear system of equations
are then approximated using the PETSc \cite{petsc-user-ref,
  petsc-efficient} Newton line search method with a tolerance of
$10^{-13}$ on each time step. A combination of Paraview and Matplotlib
have been used as visualisation tools. For each benchmark test we fix
the polynomial degree $q$ and compute a sequence of solutions with $h
= h(i) = 2^{-i}$ and $\tau = Ch$ so temporal discretisation error is
negligible.  This is done for a sequence of refinement levels, $i=l,
\dots, L$. We have previously used $S^1$ as the unitary periodic
domain. For our numerical experiments, we have scaled the domain to
$[0,40]$ for computational convenience.

\begin{Definition}[Experimental order of convergence]
  Given two sequences $a(i)$ and $h(i) \searrow 0$ we define the
  \emph{experimental order of convergence} (EOC) to be the local slope
  of the $\log{a(i)}$ vs. $\log{h(i)}$ curve, i.e.,
  \begin{equation}
    \EOC(a,h;i) = \frac{\log{\frac{a_{i+1}}{a_i}}}{\log{\frac{h_{i+1}}{h_i}}}.
  \end{equation}
\end{Definition}

\begin{Definition}[Effectivity Index] \label{def:ei}
  \revise{
  Given two sequences $a(i)$ and $b(i)$, the effectivity index is
  defined by ratio of the two, i.e.,
  \begin{equation}
    \operatorname{EI}(a,b;i)
    =
    \frac{b(i)}{a(i)}
    .
  \end{equation}
  In the sequel, we shall exclusively use the effectivity index where
  $b(i)$ is a sequence of a posteriori errors and $a(i)$ the error
  measured in the norm $\enorm{\cdot}$.  }
\end{Definition}

\begin{Remark}[Numerical deviation in $F_1$]
  While the analysis shows that our scheme \emph{exactly} preserves
  the energy over arbitrarily long time, the implementation relies on
  linear and nonlinear systems that inherently require further
  approximation. The result of this is that the energy may deviate
  \emph{locally} up to the tolerance of the linear and nonlinear
  solvers which introduces the possibility of these errors propagating
  over time. In our numerical tests we focus on studying the global
  deviation in time, $F_1(U^n)-F_1(U^0)$, which \emph{includes} any
  propagation arising from solver or precision errors.
\end{Remark}

\subsection*{Test 1: Conservativity and convergence of the linear scheme}

We begin by examining the global deviation in invariants for the
linear problem, i.e., when $f(u) = \frac12 u^2$. We observe, in Figure
\ref{fig:ldevdefocus}, that the both problems conserve the expected
invariants. 

\begin{figure}[h!]
    \subfigure[ $q=1$ ]{
\begin{tikzpicture}
  \begin{axis}[
      cycle list/Dark2,
      thick,
      xmode=linear,
      ymode=log,
      xlabel=$t_n$,
      ylabel=$\norm{F^n_{i,h} - F^0_{i,h}}$,
      ymin=1e-17, ymax=1e-3,
      grid=both,
      minor grid style={gray!25},
      major grid style={gray!25},
      width=0.45\linewidth,
      legend columns = 3,
      legend style={at={(0.725,1.05)},anchor=south west},
    ]
    \addplot
    table[x=t,y=mass, col sep=comma]{csv/dev/deviation_scheme11_ic1_degree1_coarse.csv};
    \addlegendentry{\tiny{$i=-1$}};
    \addplot
    table[x=t,y=momentum, col sep=comma]{csv/dev/deviation_scheme11_ic1_degree1_coarse.csv};
    \addlegendentry{\tiny{$i=0$}};
    \addplot
    table[x=t,y=energy, col sep=comma]{csv/dev/deviation_scheme11_ic1_degree1_coarse.csv};
    \addlegendentry{\tiny{$i=1$}};
  \end{axis}
\end{tikzpicture}
}
\hspace*{-2cm}
\hfill
\hspace*{-1cm}
  \subfigure[ $q=2$ ]{
\begin{tikzpicture}
  \begin{axis}[
      cycle list/Dark2,
      thick,
      xmode=linear,
      ymode=log,
      xlabel=$t_n$,
      ylabel=$\norm{F^n_{i,h} - F^0_{i,h}}$,
      ylabel near ticks,
      yticklabel pos=right,
      ymin=1e-17, ymax=1e-3,
      grid=both,
      minor grid style={gray!25},
      major grid style={gray!25},
      width=0.45\linewidth,
    ]
    \addplot
    table[x=t,y=mass, col sep=comma]{csv/dev/deviation_scheme11_ic1_degree2_coarse.csv};
    \addplot
    table[x=t,y=momentum, col sep=comma]{csv/dev/deviation_scheme11_ic1_degree2_coarse.csv};
    \addplot
    table[x=t,y=energy, col sep=comma]{csv/dev/deviation_scheme11_ic1_degree2_coarse.csv};
  \end{axis}
\end{tikzpicture}
}
  
  \caption{ The deviation in mass, momentum and energy with $T=100$
    for the scheme \eqref{eq:dis-mixed-system-fd} with $f(u) = \frac12
    u^2$. The initial conditions are given by
    \eqref{eqn:ldefocusexact} where $l=1,C_1=1,C_2=0$, further we
    choose $\dt{n}=0.2$, $\dx{m}=0.5$ and vary the polynomial degree
    $q$. Notice that the deviation in mass and energy is below our
    solver tolerance of $10^{-13}$ for all $q$ with the deviation in
    momentum decreasing as we increase $q$.
    \label{fig:ldevdefocus} }
\end{figure}

We plot the experimental order of convergence for the linear problem
in Figure \ref{fig:eocldefocus}. We observe that the method
convergences at the rate shown in the a priori bound
\eqref{eqn:defocusapriori}, and the a posteriori error bound
\eqref{eq:aposteriori-linear} behaves optimally.

\begin{figure}[h!]
  \centering

\subfigure[ $e$ ]{
\begin{tikzpicture}
  \begin{axis}[
      cycle list/Dark2,
      thick,
      xmode=linear,
      ymode=log,
      yscale=2,
      xlabel=$t_n$,
      ymin=1e-5, ymax=1e2,
      grid=both,
      minor grid style={gray!25},
      major grid style={gray!25},
      width=0.2\linewidth,
      legend style={at={(0,0.5)},anchor=west}]
    \addplot
    table[x=Time,y=l1, col sep=comma]{csv/con/lkdv_degree1_error_downsample.csv};
    \addplot
    table[x=Time,y=l2, col sep=comma]{csv/con/lkdv_degree1_error_downsample.csv};
    \addplot
    table[x=Time,y=l3, col sep=comma]{csv/con/lkdv_degree1_error_downsample.csv};
    \addplot
    table[x=Time,y=l4, col sep=comma]{csv/con/lkdv_degree1_error_downsample.csv};
    \addplot
    table[x=Time,y=l5, col sep=comma]{csv/con/lkdv_degree1_error_downsample.csv};
    \addplot
    table[x=Time,y=l6, col sep=comma]{csv/con/lkdv_degree1_error_downsample.csv};
  \end{axis}
\end{tikzpicture}
}
\subfigure[ $EOC(e)$ ]{
\begin{tikzpicture}
  \begin{axis}[
      cycle list/Dark2,
      thick,
      xmode=linear,
      ymode=linear,
      yscale=2,
      xlabel=$t_n$,
      ymin=0, ymax=3,
      grid=both,
      minor grid style={gray!25},
      major grid style={gray!25},
      width=0.2\linewidth,
      legend style={at={(0,0.5)},anchor=west},
    ]
    \addplot
    table[x=Time,y=l1, col sep=comma]{csv/con/lkdv_degree1_error_eoc_downsample.csv};
    \addplot
    table[x=Time,y=l2, col sep=comma]{csv/con/lkdv_degree1_error_eoc_downsample.csv};
    \addplot
    table[x=Time,y=l3, col sep=comma]{csv/con/lkdv_degree1_error_eoc_downsample.csv};
    \addplot
    table[x=Time,y=l4, col sep=comma]{csv/con/lkdv_degree1_error_eoc_downsample.csv};
    \addplot
    table[x=Time,y=l5, col sep=comma]{csv/con/lkdv_degree1_error_eoc_downsample.csv};
    \addplot
    table[x=Time,y=l6, col sep=comma]{csv/con/lkdv_degree1_error_eoc_downsample.csv};
  \end{axis}
\end{tikzpicture}
}
\hspace*{-0.5cm}
\subfigure[ $E_{st}$ ]{
\begin{tikzpicture}
  \begin{axis}[
      cycle list/Dark2,
      thick,
      xmode=linear,
      ymode=log,
      yscale=2,
      xlabel=$t_n$,
      ymin=1e-5, ymax=1e2,
      grid=both,
      minor grid style={gray!25},
      major grid style={gray!25},
      width=0.2\linewidth,
      legend columns = 3,
      legend style={at={(-0.6,1)},anchor=north west,nodes={scale=0.75, transform shape}},
    ]
    \addplot
    table[x=Time,y=l1, col sep=comma]{csv/con/lkdv_degree1_H_downsample.csv};
    \addlegendentry{\scriptsize{$h_1=1.25$}};
    \addplot
    table[x=Time,y=l2, col
    sep=comma]{csv/con/lkdv_degree1_H_downsample.csv};
    \addlegendentry{\scriptsize{$h_2=0.63$}};
    \addplot
    table[x=Time,y=l3, col
    sep=comma]{csv/con/lkdv_degree1_H_downsample.csv};
    \addlegendentry{\scriptsize{$h_3=0.31$}};
    \addplot
    table[x=Time,y=l4, col
    sep=comma]{csv/con/lkdv_degree1_H_downsample.csv};
    \addlegendentry{\scriptsize{$h_4=0.16$}};
    \addplot
    table[x=Time,y=l5, col
    sep=comma]{csv/con/lkdv_degree1_H_downsample.csv};
    \addlegendentry{\scriptsize{$h_5=0.08$}};
    \addplot
    table[x=Time,y=l6, col
    sep=comma]{csv/con/lkdv_degree1_H_downsample.csv};
    \addlegendentry{\scriptsize{$h_6=0.04$}};
  \end{axis}
\end{tikzpicture}
}
\hspace*{-1.5cm}
\subfigure[ $EOC(E_{st})$ ]{
\begin{tikzpicture}
  \begin{axis}[
      cycle list/Dark2,
      thick,
      xmode=linear,
      ymode=linear,
      yscale=2,
      xlabel=$t_n$,
      ymin=0, ymax=3,
      grid=both,
      minor grid style={gray!25},
      major grid style={gray!25},
      width=0.2\linewidth,
      legend style={at={(0,0.5)},anchor=west},
    ]
    \addplot
    table[x=Time,y=l1, col sep=comma]{csv/con/lkdv_degree1_H_eoc_downsample.csv};
    \addplot
    table[x=Time,y=l2, col sep=comma]{csv/con/lkdv_degree1_H_eoc_downsample.csv};
    \addplot
    table[x=Time,y=l3, col sep=comma]{csv/con/lkdv_degree1_H_eoc_downsample.csv};
    \addplot
    table[x=Time,y=l4, col sep=comma]{csv/con/lkdv_degree1_H_eoc_downsample.csv};
    \addplot
    table[x=Time,y=l5, col sep=comma]{csv/con/lkdv_degree1_H_eoc_downsample.csv};
    \addplot
    table[x=Time,y=l6, col sep=comma]{csv/con/lkdv_degree1_H_eoc_downsample.csv};
  \end{axis}
\end{tikzpicture}
}
\subfigure[ $EI$ ]{
\begin{tikzpicture}
  \begin{axis}[
      cycle list/Dark2,
      thick,
      xmode=linear,
      ymode=linear,
      yscale=2,
      xlabel=$t_n$,
      grid=both,
      minor grid style={gray!25},
      major grid style={gray!25},
      width=0.2\linewidth,
      legend style={at={(0,0.5)},anchor=west},
    ]
    \addplot
    table[x=Time,y=l1, col sep=comma]{csv/con/lkdv_degree1_EI_downsample.csv};
    \addplot
    table[x=Time,y=l2, col sep=comma]{csv/con/lkdv_degree1_EI_downsample.csv};
    \addplot
    table[x=Time,y=l3, col sep=comma]{csv/con/lkdv_degree1_EI_downsample.csv};
    \addplot
    table[x=Time,y=l4, col sep=comma]{csv/con/lkdv_degree1_EI_downsample.csv};
    \addplot
    table[x=Time,y=l5, col sep=comma]{csv/con/lkdv_degree1_EI_downsample.csv};
    \addplot
    table[x=Time,y=l6, col sep=comma]{csv/con/lkdv_degree1_EI_downsample.csv};
  \end{axis}
\end{tikzpicture}
}
  \caption{The error $e:=\enorm{u-U}$ and the a posteriori estimator
    (\ref{eq:aposteriori-linear2}) denoted $E_{st}$, and their
    associated experimental order of convergence with the
    corresponding exact solution \eqref{eqn:ldefocusexact} where
    $l=1, C_1=1, C_2=0$ with polynomial degree $q=1$. We vary
    $\dx{m} = 0.1 \dt{n}$, and observe the a priori bound
    \eqref{eqn:defocusapriori} is attained for $q=1$, and that the a
    posteriori bound is optimal. \label{fig:eocldefocus} }
\end{figure}

\begin{figure}[h!]
  \centering

\subfigure[ $e$ ]{
\begin{tikzpicture}
  \begin{axis}[
      cycle list/Dark2,
      thick,
      xmode=linear,
      ymode=log,
      yscale=2,
      xlabel=$t_n$,
      ymin=1e-5, ymax=1e2,
      grid=both,
      minor grid style={gray!25},
      major grid style={gray!25},
      width=0.2\linewidth,
      legend style={at={(0,0.5)},anchor=west}]
    \addplot
    table[x=Time,y=l1, col sep=comma]{csv/con/lkdv_degree2_error_downsample.csv};
    \addplot
    table[x=Time,y=l2, col sep=comma]{csv/con/lkdv_degree2_error_downsample.csv};
    \addplot
    table[x=Time,y=l3, col sep=comma]{csv/con/lkdv_degree2_error_downsample.csv};
    \addplot
    table[x=Time,y=l4, col sep=comma]{csv/con/lkdv_degree2_error_downsample.csv};
    \addplot
    table[x=Time,y=l5, col sep=comma]{csv/con/lkdv_degree2_error_downsample.csv};
    \addplot
    table[x=Time,y=l6, col sep=comma]{csv/con/lkdv_degree2_error_downsample.csv};
  \end{axis}
\end{tikzpicture}
}
\subfigure[ $EOC(e)$ ]{
\begin{tikzpicture}
  \begin{axis}[
      cycle list/Dark2,
      thick,
      xmode=linear,
      ymode=linear,
      yscale=2,
      xlabel=$t_n$,
      ymin=0, ymax=3,
      grid=both,
      minor grid style={gray!25},
      major grid style={gray!25},
      width=0.2\linewidth,
      legend style={at={(0,0.5)},anchor=west},
    ]
    \addplot
    table[x=Time,y=l1, col sep=comma]{csv/con/lkdv_degree2_error_eoc_downsample.csv};
    \addplot
    table[x=Time,y=l2, col sep=comma]{csv/con/lkdv_degree2_error_eoc_downsample.csv};
    \addplot
    table[x=Time,y=l3, col sep=comma]{csv/con/lkdv_degree2_error_eoc_downsample.csv};
    \addplot
    table[x=Time,y=l4, col sep=comma]{csv/con/lkdv_degree2_error_eoc_downsample.csv};
    \addplot
    table[x=Time,y=l5, col sep=comma]{csv/con/lkdv_degree2_error_eoc_downsample.csv};
    \addplot
    table[x=Time,y=l6, col sep=comma]{csv/con/lkdv_degree2_error_eoc_downsample.csv};
  \end{axis}
\end{tikzpicture}
}
\hspace*{-0.5cm}
\subfigure[ $E_{st}$ ]{
\begin{tikzpicture}
  \begin{axis}[
      cycle list/Dark2,
      thick,
      xmode=linear,
      ymode=log,
      yscale=2,
      xlabel=$t_n$,
      ymin=1e-5, ymax=1e2,
      grid=both,
      minor grid style={gray!25},
      major grid style={gray!25},
      width=0.2\linewidth,
      legend columns = 3,
      legend style={at={(-0.6,1)},anchor=north west,nodes={scale=0.75, transform shape}},
    ]
    \addplot
    table[x=Time,y=l1, col sep=comma]{csv/con/lkdv_degree2_H_downsample.csv};
    \addlegendentry{\scriptsize{$h_1=1.25$}};
    \addplot
    table[x=Time,y=l2, col
    sep=comma]{csv/con/lkdv_degree2_H_downsample.csv};
    \addlegendentry{\scriptsize{$h_2=0.63$}};
    \addplot
    table[x=Time,y=l3, col
    sep=comma]{csv/con/lkdv_degree2_H_downsample.csv};
    \addlegendentry{\scriptsize{$h_3=0.31$}};
    \addplot
    table[x=Time,y=l4, col
    sep=comma]{csv/con/lkdv_degree2_H_downsample.csv};
    \addlegendentry{\scriptsize{$h_4=0.16$}};
    \addplot
    table[x=Time,y=l5, col
    sep=comma]{csv/con/lkdv_degree2_H_downsample.csv};
    \addlegendentry{\scriptsize{$h_5=0.08$}};
    \addplot
    table[x=Time,y=l6, col
    sep=comma]{csv/con/lkdv_degree2_H_downsample.csv};
    \addlegendentry{\scriptsize{$h_6=0.04$}};
  \end{axis}
\end{tikzpicture}
}
\hspace*{-1.5cm}
\subfigure[ $EOC(E_{st})$ ]{
\begin{tikzpicture}
  \begin{axis}[
      cycle list/Dark2,
      thick,
      xmode=linear,
      ymode=linear,
      yscale=2,
      xlabel=$t_n$,
      ymin=0, ymax=3,
      grid=both,
      minor grid style={gray!25},
      major grid style={gray!25},
      width=0.2\linewidth,
      legend style={at={(0,0.5)},anchor=west},
    ]
    \addplot
    table[x=Time,y=l1, col sep=comma]{csv/con/lkdv_degree2_H_eoc_downsample.csv};
    \addplot
    table[x=Time,y=l2, col sep=comma]{csv/con/lkdv_degree2_H_eoc_downsample.csv};
    \addplot
    table[x=Time,y=l3, col sep=comma]{csv/con/lkdv_degree2_H_eoc_downsample.csv};
    \addplot
    table[x=Time,y=l4, col sep=comma]{csv/con/lkdv_degree2_H_eoc_downsample.csv};
    \addplot
    table[x=Time,y=l5, col sep=comma]{csv/con/lkdv_degree2_H_eoc_downsample.csv};
    \addplot
    table[x=Time,y=l6, col sep=comma]{csv/con/lkdv_degree2_H_eoc_downsample.csv};
  \end{axis}
\end{tikzpicture}
}
\subfigure[ $EI$ ]{
\begin{tikzpicture}
  \begin{axis}[
      cycle list/Dark2,
      thick,
      xmode=linear,
      ymode=linear,
      yscale=2,
      xlabel=$t_n$,
      grid=both,
      minor grid style={gray!25},
      major grid style={gray!25},
      width=0.2\linewidth,
      legend style={at={(0,0.5)},anchor=west},
    ]
    \addplot
    table[x=Time,y=l1, col sep=comma]{csv/con/lkdv_degree2_EI_downsample.csv};
    \addplot
    table[x=Time,y=l2, col sep=comma]{csv/con/lkdv_degree2_EI_downsample.csv};
    \addplot
    table[x=Time,y=l3, col sep=comma]{csv/con/lkdv_degree2_EI_downsample.csv};
    \addplot
    table[x=Time,y=l4, col sep=comma]{csv/con/lkdv_degree2_EI_downsample.csv};
    \addplot
    table[x=Time,y=l5, col sep=comma]{csv/con/lkdv_degree2_EI_downsample.csv};
    \addplot
    table[x=Time,y=l6, col sep=comma]{csv/con/lkdv_degree2_EI_downsample.csv};
  \end{axis}
\end{tikzpicture}
}
  \caption{The error $e:=\enorm{u-U}$ and the a posteriori estimator
    \eqref{eq:aposteriori-linear2} denoted $E_{st}$, and their
    associated experimental order of convergence with the
    corresponding exact solution \eqref{eqn:ldefocusexact} where
    $l=1, C_1=1, C_2=0$ with polynomial degree $q=2$. We vary
    $\dx{m} = 0.1 \dt{n}$, and observe the a priori bound
    \eqref{eqn:defocusapriori} is attained for $q=1$, and that the a
    posteriori bound is optimal. \label{fig:eocldefocus} }
\end{figure}

\subsection*{Test 2: Conservativity and convergence of the nonlinear scheme}

Through examining the global deviation in invariants for the nonlinear
problem $f(u) = \frac12 u^4$ we observe, in Figure
\ref{fig:devdefocus}, that the both problems conserve the expected
invariants. 

\begin{figure}[h!]
    \subfigure[ $q=1$ ]{
\begin{tikzpicture}
  \begin{axis}[
      cycle list/Dark2,
      thick,
      xmode=linear,
      ymode=log,
      xlabel=$t_n$,
      ylabel=$\norm{F^n_{i,h} - F^0_{i,h}}$,
      ymin=1e-17, ymax=1e-2,
      grid=both,
      minor grid style={gray!25},
      major grid style={gray!25},
      width=0.45\linewidth,
      legend columns = 3,
      legend style={at={(0.725,1.05)},anchor=south west},
    ]
    \addplot
    table[x=t,y=mass, col sep=comma]{csv/dev/deviation_scheme12_ic1_degree1_coarse.csv};
    \addlegendentry{\tiny{$i=-1$}};
    \addplot
    table[x=t,y=momentum, col sep=comma]{csv/dev/deviation_scheme12_ic1_degree1_coarse.csv};
    \addlegendentry{\tiny{$i=0$}};
    \addplot
    table[x=t,y=energy, col sep=comma]{csv/dev/deviation_scheme12_ic1_degree1_coarse.csv};
    \addlegendentry{\tiny{$i=1$}};
  \end{axis}
\end{tikzpicture}
  }
\hspace*{-2cm}
\hfill
\hspace*{-1cm}
  \subfigure[ $q=2$ ]{
\begin{tikzpicture}
  \begin{axis}[
      cycle list/Dark2,
      thick,
      xmode=linear,
      ymode=log,
      xlabel=$t_n$,
      ylabel=$\norm{F^n_{i,h} - F^0_{i,h}}$,
      ylabel near ticks,
      yticklabel pos=right,
      ymin=1e-17, ymax=1e-2,
      grid=both,
      minor grid style={gray!25},
      major grid style={gray!25},
      width=0.45\linewidth,
    ]
    \addplot
    table[x=t,y=mass, col sep=comma]{csv/dev/deviation_scheme12_ic1_degree2_coarse.csv};
    \addplot
    table[x=t,y=momentum, col sep=comma]{csv/dev/deviation_scheme12_ic1_degree2_coarse.csv};
    \addplot
    table[x=t,y=energy, col sep=comma]{csv/dev/deviation_scheme12_ic1_degree2_coarse.csv};
  \end{axis}
\end{tikzpicture}
}
  
  \caption{ The deviation in mass, momentum and energy with $T=100$
    for the scheme \eqref{eq:dis-mixed-system-fd} with $f(u) = \frac12
    u^4$. The initial conditions are given by \eqref{eqn:defocusexact}
    with $k=0.9$. We additionally stretch our spatial interval to $x
    \in [0,16 K(k)]$, where $K(k)$ is the complete elliptic integral
    of the first kind, which numerically ensures that our solution is
    periodic up to a tolerance of $10^{-13}$.  Further we choose
    $\dt{n}=0.2$, $\dx{m}=0.5$ and vary the polynomial degree $q$. The
    deviation in both mass and energy is below solver precision, with
    momentum decreasing as we increase $q$. We note that as we
    increase $q$ the momentum does not decrease as quickly as for the
    linear case.
    \label{fig:devdefocus} }
\end{figure}

In addition, we benchmark the nonlinear scheme against the exact
solution \eqref{eqn:defocusexact} with $k=0.9$ over the stretched
spatial domain $x \in [0,16 K(k)]$ (where $K(k)$ is the complete
elliptic integral of the first kind) yielding Figure
\ref{fig:eocdefocus}. We observe similar convergence rates to the
linear case, satisfying the a posteriori error bound, in addition to
indicating the existence of optimal a priori bounds.

\begin{figure}[h!]
  \centering

\subfigure[ $e$ ]{
\begin{tikzpicture}
  \begin{axis}[
      cycle list/Dark2,
      thick,
      xmode=linear,
      ymode=log,
      yscale=2,
      xlabel=$t_n$,
      ymin=1e-6, ymax=1e2,
      grid=both,
      minor grid style={gray!25},
      major grid style={gray!25},
      width=0.2\linewidth,
      legend style={at={(0,0.5)},anchor=west}]
    \addplot
    table[x=Time,y=l1, col sep=comma]{csv/con/defocus_degree1_error_downsample.csv};
    \addplot
    table[x=Time,y=l2, col sep=comma]{csv/con/defocus_degree1_error_downsample.csv};
    \addplot
    table[x=Time,y=l3, col sep=comma]{csv/con/defocus_degree1_error_downsample.csv};
    \addplot
    table[x=Time,y=l4, col sep=comma]{csv/con/defocus_degree1_error_downsample.csv};
    \addplot
    table[x=Time,y=l5, col sep=comma]{csv/con/defocus_degree1_error_downsample.csv};
    \addplot
    table[x=Time,y=l6, col sep=comma]{csv/con/defocus_degree1_error_downsample.csv};
  \end{axis}
\end{tikzpicture}
}
\subfigure[ $EOC(e)$ ]{
\begin{tikzpicture}
  \begin{axis}[
      cycle list/Dark2,
      thick,
      xmode=linear,
      ymode=linear,
      yscale=2,
      xlabel=$t_n$,
      ymin=0, ymax=3,
      grid=both,
      minor grid style={gray!25},
      major grid style={gray!25},
      width=0.2\linewidth,
      legend style={at={(0,0.5)},anchor=west},
    ]
    \addplot
    table[x=Time,y=l1, col sep=comma]{csv/con/defocus_degree1_error_eoc_downsample.csv};
    \addplot
    table[x=Time,y=l2, col sep=comma]{csv/con/defocus_degree1_error_eoc_downsample.csv};
    \addplot
    table[x=Time,y=l3, col sep=comma]{csv/con/defocus_degree1_error_eoc_downsample.csv};
    \addplot
    table[x=Time,y=l4, col sep=comma]{csv/con/defocus_degree1_error_eoc_downsample.csv};
    \addplot
    table[x=Time,y=l5, col sep=comma]{csv/con/defocus_degree1_error_eoc_downsample.csv};
    \addplot
    table[x=Time,y=l6, col sep=comma]{csv/con/defocus_degree1_error_eoc_downsample.csv};
  \end{axis}
\end{tikzpicture}
}
\hspace*{-0.5cm}
\subfigure[ $E_{st}$ ]{
\begin{tikzpicture}
  \begin{axis}[
      cycle list/Dark2,
      thick,
      xmode=linear,
      ymode=log,
      yscale=2,
      xlabel=$t_n$,
      ymin=1e-6, ymax=1e2,
      grid=both,
      minor grid style={gray!25},
      major grid style={gray!25},
      width=0.2\linewidth,
      legend columns = 3,
      legend style={at={(-0.6,1.02)},anchor=north west,nodes={scale=0.75, transform shape}},
    ]
    \addplot
    table[x=Time,y=l1, col sep=comma]{csv/con/defocus_degree1_H_downsample.csv};
    \addlegendentry{\scriptsize{$h_1=0.31$}};
    \addplot
    table[x=Time,y=l2, col
    sep=comma]{csv/con/defocus_degree1_H_downsample.csv};
    \addlegendentry{\scriptsize{$h_2=0.16$}};
    \addplot
    table[x=Time,y=l3, col
    sep=comma]{csv/con/defocus_degree1_H_downsample.csv};
    \addlegendentry{\scriptsize{$h_3=0.08$}};
    \addplot
    table[x=Time,y=l4, col
    sep=comma]{csv/con/defocus_degree1_H_downsample.csv};
    \addlegendentry{\scriptsize{$h_4=0.04$}};
    \addplot
    table[x=Time,y=l5, col
    sep=comma]{csv/con/defocus_degree1_H_downsample.csv};
    \addlegendentry{\scriptsize{$h_5=0.02$}};
    \addplot
    table[x=Time,y=l6, col
    sep=comma]{csv/con/defocus_degree1_H_downsample.csv};
    \addlegendentry{\scriptsize{$h_6=0.01$}};
  \end{axis}
\end{tikzpicture}
}
\hspace*{-1.5cm}
\subfigure[ $EOC(E_{st})$ ]{
\begin{tikzpicture}
  \begin{axis}[
      cycle list/Dark2,
      thick,
      xmode=linear,
      ymode=linear,
      yscale=2,
      xlabel=$t_n$,
      ymin=0, ymax=3,
      grid=both,
      minor grid style={gray!25},
      major grid style={gray!25},
      width=0.2\linewidth,
      legend style={at={(0,0.5)},anchor=west},
    ]
    \addplot
    table[x=Time,y=l1, col sep=comma]{csv/con/defocus_degree1_H_eoc_downsample.csv};
    \addplot
    table[x=Time,y=l2, col sep=comma]{csv/con/defocus_degree1_H_eoc_downsample.csv};
    \addplot
    table[x=Time,y=l3, col sep=comma]{csv/con/defocus_degree1_H_eoc_downsample.csv};
    \addplot
    table[x=Time,y=l4, col sep=comma]{csv/con/defocus_degree1_H_eoc_downsample.csv};
    \addplot
    table[x=Time,y=l5, col sep=comma]{csv/con/defocus_degree1_H_eoc_downsample.csv};
    \addplot
    table[x=Time,y=l6, col sep=comma]{csv/con/defocus_degree1_H_eoc_downsample.csv};
  \end{axis}
\end{tikzpicture}
}
\subfigure[ $EI$ ]{
\begin{tikzpicture}
  \begin{axis}[
      cycle list/Dark2,
      thick,
      xmode=linear,
      ymode=linear,
      yscale=2,
      xlabel=$t_n$,
      grid=both,
      minor grid style={gray!25},
      major grid style={gray!25},
      width=0.2\linewidth,
      legend style={at={(0,0.5)},anchor=west},
    ]
    \addplot
    table[x=Time,y=l1, col sep=comma]{csv/con/defocus_degree1_EI_downsample.csv};
    \addplot
    table[x=Time,y=l2, col sep=comma]{csv/con/defocus_degree1_EI_downsample.csv};
    \addplot
    table[x=Time,y=l3, col sep=comma]{csv/con/defocus_degree1_EI_downsample.csv};
    \addplot
    table[x=Time,y=l4, col sep=comma]{csv/con/defocus_degree1_EI_downsample.csv};
    \addplot
    table[x=Time,y=l5, col sep=comma]{csv/con/defocus_degree1_EI_downsample.csv};
    \addplot
    table[x=Time,y=l6, col sep=comma]{csv/con/defocus_degree1_EI_downsample.csv};
  \end{axis}
\end{tikzpicture}
}
  \caption{ The error $e:=\enorm{u-U}$ and the a posteriori estimator
    \eqref{eq:aposteriori-nonlinear} denoted $E_{st}$, and their
    associated experimental order of convergence using the solution
    \eqref{eqn:defocusexact} with $k=0.9$ and using polynomial degree
    $q=1$. We vary $\dx{m} = 0.1 \dt{n}$, and observe that the a
    priori and a posteriori error bounds both converge
    optimally. \label{fig:eocdefocus} }
\end{figure}

\begin{figure}[h!]
  \centering

\subfigure[ $e$ ]{
\begin{tikzpicture}
  \begin{axis}[
      cycle list/Dark2,
      thick,
      xmode=linear,
      ymode=log,
      yscale=2,
      xlabel=$t_n$,
      ymin=1e-6, ymax=1e2,
      grid=both,
      minor grid style={gray!25},
      major grid style={gray!25},
      width=0.2\linewidth,
      legend style={at={(0,0.5)},anchor=west}]
    \addplot
    table[x=Time,y=l1, col sep=comma]{csv/con/defocus_degree2_error_downsample.csv};
    \addplot
    table[x=Time,y=l2, col sep=comma]{csv/con/defocus_degree2_error_downsample.csv};
    \addplot
    table[x=Time,y=l3, col sep=comma]{csv/con/defocus_degree2_error_downsample.csv};
    \addplot
    table[x=Time,y=l4, col sep=comma]{csv/con/defocus_degree2_error_downsample.csv};
    \addplot
    table[x=Time,y=l5, col sep=comma]{csv/con/defocus_degree2_error_downsample.csv};
    \addplot
    table[x=Time,y=l6, col sep=comma]{csv/con/defocus_degree2_error_downsample.csv};
  \end{axis}
\end{tikzpicture}
}
\subfigure[ $EOC(e)$ ]{
\begin{tikzpicture}
  \begin{axis}[
      cycle list/Dark2,
      thick,
      xmode=linear,
      ymode=linear,
      yscale=2,
      xlabel=$t_n$,
      ymin=0, ymax=3,
      grid=both,
      minor grid style={gray!25},
      major grid style={gray!25},
      width=0.2\linewidth,
      legend style={at={(0,0.5)},anchor=west},
    ]
    \addplot
    table[x=Time,y=l1, col sep=comma]{csv/con/defocus_degree2_error_eoc_downsample.csv};
    \addplot
    table[x=Time,y=l2, col sep=comma]{csv/con/defocus_degree2_error_eoc_downsample.csv};
    \addplot
    table[x=Time,y=l3, col sep=comma]{csv/con/defocus_degree2_error_eoc_downsample.csv};
    \addplot
    table[x=Time,y=l4, col sep=comma]{csv/con/defocus_degree2_error_eoc_downsample.csv};
    \addplot
    table[x=Time,y=l5, col sep=comma]{csv/con/defocus_degree2_error_eoc_downsample.csv};
    \addplot
    table[x=Time,y=l6, col sep=comma]{csv/con/defocus_degree2_error_eoc_downsample.csv};
  \end{axis}
\end{tikzpicture}
}
\hspace*{-0.5cm}
\subfigure[ $E_{st}$ ]{
\begin{tikzpicture}
  \begin{axis}[
      cycle list/Dark2,
      thick,
      xmode=linear,
      ymode=log,
      yscale=2,
      xlabel=$t_n$,
      ymin=1e-6, ymax=1e2,
      grid=both,
      minor grid style={gray!25},
      major grid style={gray!25},
      width=0.2\linewidth,
      legend columns = 3,
      legend style={at={(-0.6,1.02)},anchor=north west,nodes={scale=0.75, transform shape}},
    ]
    \addplot
    table[x=Time,y=l1, col sep=comma]{csv/con/defocus_degree2_H_downsample.csv};
    \addlegendentry{\scriptsize{$h_1=0.31$}};
    \addplot
    table[x=Time,y=l2, col
    sep=comma]{csv/con/defocus_degree2_H_downsample.csv};
    \addlegendentry{\scriptsize{$h_2=0.16$}};
    \addplot
    table[x=Time,y=l3, col
    sep=comma]{csv/con/defocus_degree2_H_downsample.csv};
    \addlegendentry{\scriptsize{$h_3=0.08$}};
    \addplot
    table[x=Time,y=l4, col
    sep=comma]{csv/con/defocus_degree2_H_downsample.csv};
    \addlegendentry{\scriptsize{$h_4=0.04$}};
    \addplot
    table[x=Time,y=l5, col
    sep=comma]{csv/con/defocus_degree2_H_downsample.csv};
    \addlegendentry{\scriptsize{$h_5=0.02$}};
    \addplot
    table[x=Time,y=l6, col
    sep=comma]{csv/con/defocus_degree2_H_downsample.csv};
    \addlegendentry{\scriptsize{$h_6=0.01$}};
  \end{axis}
\end{tikzpicture}
}
\hspace*{-1.5cm}
\subfigure[ $EOC(E_{st})$ ]{
\begin{tikzpicture}
  \begin{axis}[
      cycle list/Dark2,
      thick,
      xmode=linear,
      ymode=linear,
      yscale=2,
      xlabel=$t_n$,
      ymin=0, ymax=3,
      grid=both,
      minor grid style={gray!25},
      major grid style={gray!25},
      width=0.2\linewidth,
      legend style={at={(0,0.5)},anchor=west},
    ]
    \addplot
    table[x=Time,y=l1, col sep=comma]{csv/con/defocus_degree2_H_eoc_downsample.csv};
    \addplot
    table[x=Time,y=l2, col sep=comma]{csv/con/defocus_degree2_H_eoc_downsample.csv};
    \addplot
    table[x=Time,y=l3, col sep=comma]{csv/con/defocus_degree2_H_eoc_downsample.csv};
    \addplot
    table[x=Time,y=l4, col sep=comma]{csv/con/defocus_degree2_H_eoc_downsample.csv};
    \addplot
    table[x=Time,y=l5, col sep=comma]{csv/con/defocus_degree2_H_eoc_downsample.csv};
    \addplot
    table[x=Time,y=l6, col sep=comma]{csv/con/defocus_degree2_H_eoc_downsample.csv};
  \end{axis}
\end{tikzpicture}
}
\subfigure[ $EI$ ]{
\begin{tikzpicture}
  \begin{axis}[
      cycle list/Dark2,
      thick,
      xmode=linear,
      ymode=linear,
      yscale=2,
      xlabel=$t_n$,
      grid=both,
      minor grid style={gray!25},
      major grid style={gray!25},
      width=0.2\linewidth,
      legend style={at={(0,0.5)},anchor=west},
    ]
    \addplot
    table[x=Time,y=l1, col sep=comma]{csv/con/defocus_degree2_EI_downsample.csv};
    \addplot
    table[x=Time,y=l2, col sep=comma]{csv/con/defocus_degree2_EI_downsample.csv};
    \addplot
    table[x=Time,y=l3, col sep=comma]{csv/con/defocus_degree2_EI_downsample.csv};
    \addplot
    table[x=Time,y=l4, col sep=comma]{csv/con/defocus_degree2_EI_downsample.csv};
    \addplot
    table[x=Time,y=l5, col sep=comma]{csv/con/defocus_degree2_EI_downsample.csv};
    \addplot
    table[x=Time,y=l6, col sep=comma]{csv/con/defocus_degree2_EI_downsample.csv};
  \end{axis}
\end{tikzpicture}
}
  \caption{ The error $e:=\enorm{u-U}$ and the a posteriori estimator
    \eqref{eq:aposteriori-nonlinear} denoted $E_{st}$, and their
    associated experimental order of convergence using the solution
    \eqref{eqn:defocusexact} with $k=0.9$ and using polynomial degree
    $q=2$. We vary $\dx{m} = 0.1 \dt{n}$, and observe that the a
    priori and a posteriori error bounds both converge
    optimally. \label{fig:eocdefocus} }
\end{figure}

\section{Conclusions and outlook}

In this work we have developed a novel discontinuous Galerkin scheme
for a specific class of Hamiltonian problem. We have shown the method
to be well posed and demonstrated that it inherits desirable
conservative properties of the PDE.

Furthermore, we have conducted an a priori error analysis that shows
the method is optimally convergent in the energy norm. This is not
surprising as the method is deliberately designed to be conservative
and the invariant induces the energy norm.

An a posteriori analysis was carried out for the semi discrete scheme
that is very much in the spirit of the original elliptic
reconstruction framework of \cite{Makridakis:2003}. It is shown that
we can make use of this framework to derive a posteriori bounds in the
energy norm, this is different to the framework of
\cite{KarakashianMakridakis:2015} where a dispersive reconstruction
was used to enable $\leb{2}$ error control. An attractive feature of
the analysis we present is that the bound holds irrespective of the
underlying polynomial degree of the approximation scheme.

As an outlook we plan to make use of the a posteriori framework
developed here to extend to fully discrete a posteriori bounds that
are able to account for mesh adaptivity. This is a particularly subtle
point as mesh change, when done in a naive way, can actually induce
instabilities \cite{BanschKarakatsaniMakridakis:2013} although one can
design adaptive schemes that ensure compatibility with the underlying
Hamiltonian formulation of the problem \cite{Eidnes:2018,
  MiyatakeMatsuo:2015}.

\clearpage

\bibliographystyle{abbrv}      
\bibliography{central,tristansbib,tristanswritings}   

\end{document}